\documentclass[twocolumn]{autart}
\usepackage{bm}
\usepackage{graphics}
\usepackage{epsfig}
\usepackage{times}
\usepackage{amsmath}
\usepackage{amssymb}
\usepackage{graphicx}
\usepackage{epstopdf}
\usepackage{epsfig}
\usepackage{enumerate}
\usepackage[noadjust]{cite}
\usepackage{color}
\usepackage{makecell}
\usepackage{arydshln}
\usepackage{stmaryrd}
\usepackage{stackengine}
\usepackage{threeparttable}
\usepackage{algorithm}
\usepackage{algorithmicx}
\usepackage{xpatch}
\usepackage{bbm}
\usepackage{bbding}
\usepackage{enumerate}
\usepackage{algpseudocode}
\usepackage{multirow}
\usepackage{epsfig}
\usepackage{tikz,harvard,fancybox}

\makeatletter
\xpatchcmd\ALG@step{\arabic{ALG@line}}{\fmtlinenumber{ALG@line}}{}{}
\makeatother
\let\fmtlinenumber\arabic %by default, line numbers are arabic numbers
\newtheorem{Theorem}{Theorem}
\newtheorem{Assumption}{Assumption}
\newtheorem{Corollary}{Corollary}
\newtheorem{Lemma}{Lemma}

\newfont{\bb}{msbm10}

\newcommand{\tr}{^{\sf T}}
\newcommand{\M}[1]{{\bf{#1}}}

\newcommand{\m}[1]{{\mathrm{#1}}}

\allowdisplaybreaks[4]
\begin{document}
\begin{frontmatter}

\title{BALPA: A Balanced Primal--Dual Algorithm for Nonsmooth Optimization with Application to Distributed Optimization\thanksref{footnoteinfo1}} % Title, preferably not more

\thanks[footnoteinfo1]{This work was supported in part by the National Natural Science Foundation of China under Grant Nos. 61833005, 62003084 and 62176056, the Natural Science Foundation of Jiangsu Province of China under Grant No. BK20200355, and Young Elite Scientists Sponsorship Program by CAST, 2021QNRC001.}

\thanks[footnoteinfo2]{Corresponding author: Jinde Cao.}

\author[Math]{Luyao Guo}\ead{ly\_guo@seu.edu.cn},    % Add the
\author[Math]{Jinde Cao\thanksref{footnoteinfo2}}\ead{jdcao@seu.edu.cn},
\author[Cyber]{Xinli Shi}\ead{xinli\_shi@seu.edu.cn},
\author[Comp]{Shaofu Yang}\ead{sfyang@seu.edu.cn}.
%\author[Cyber]{Xinli Shi}\ead{xinli\_shi@seu.edu.cn}.         % e-mail address

\address[Math]{School of Mathematics, Southeast University, Nanjing 210096, China}  % Please supply
\address[Cyber]{School of Cyber Science \& Engineering, Southeast University, Nanjing 210096, China}             % full addresses
\address[Comp]{School of Computer Science and Engineering, Southeast University, Nanjing 210096, China}
%\address[Baiae]{Yonsei Frontier Lab, Yonsei University, Seoul, Korea}        % here.

\begin{keyword}                           % Five to ten keywords,
Composite optimization, distributed optimization, primal--dual algorithm, stochastic algorithms.
% chosen from the IFAC
\end{keyword}                             % keyword list or with the
                                          % help of the Automatica
                                          % keyword wizard

\begin{abstract}                          % Abstract of not more than 200 words.
In this paper, we propose a novel primal--dual proximal splitting algorithm (PD-PSA), named BALPA, for the composite optimization problem with equality constraints, where the loss function consists of a smooth term and a nonsmooth term composed with a linear mapping. In BALPA, the dual update is designed as a proximal point for a time-varying quadratic function, which balances the implementation of primal and dual update and retains the proximity-induced feature of classic PD-PSAs. In addition, by this balance, BALPA eliminates the inefficiency of classic PD-PSAs for composite optimization problems in which the Euclidean norm of the linear mapping or the equality constraint mapping is large. Therefore, BALPA not only inherits the advantages of simple structure and easy implementation of classic PD-PSAs but also ensures a fast convergence when these norms are large. Moreover, we propose a stochastic version of BALPA (S-BALPA) and apply the developed BALPA to distributed optimization to devise a new distributed optimization algorithm. Furthermore, a comprehensive convergence analysis for BALPA and S-BALPA is conducted, respectively. Finally, numerical experiments demonstrate the efficiency of the proposed algorithms.
\end{abstract}

\end{frontmatter}

%%%%%%%%%%%%%%%%%%%%%%%%%%%%%%%%%%%%%%%%%%%%%%%%%%%%%%%%%%%%%%%%%%%%%%%%%%%%
\section{Introduction}
Consider the following composite optimization problem
\begin{equation}\label{Problem}
\min_{\m{x}} f(\m{x})+r(\m{B}\m{x})+\delta_0(\m{Dx}-\m{d}),
\end{equation}
where $\m{x}\in\mathbb{R}^n$, $\m{B}\in\mathbb{R}^{p_1\times n}$ and $\m{D}\in\mathbb{R}^{p_2\times n}$ are two given matrices, $\m{d}\in\mathbb{R}^{p_2}$, $f:\mathbb{R}^n\rightarrow(-\infty,+\infty)$ is convex, $r:\mathbb{R}^{p_1}\rightarrow(-\infty,+\infty]$ is a proper closed convex function, and $\delta_0$ is an indicator function defined as $\delta_0(\m{Dx}-\m{d})=0$ if $\m{Dx}-\m{d}=0$, otherwise $\delta_0(\m{Dx}-\m{d})=\infty$, which encodes the equality constraint $\m{Dx-d}=0$. We assume that $f$ is differentiable with $L$-Lipschitz continuous gradient and $r$ is proximable, i.e., the proximal mapping of $r$ defined as
$
\mathrm{prox}^{\alpha}_{r}(\m{y})=\arg\min_{\m{v}} r(\m{v})+\frac{1}{2\alpha}\|\m{v}-\m{y}\|^2,
$
has analytical solution or can be computed efficiently.

This composite structure covers a wide range of optimization scenarios, such as optimal voltage control in networks \cite{Auto1}, complete vehicle energy management \cite{Auto2}, and distributed optimization or learning \cite{Nedic20151}. In particular, we take the distributed learning as an example. Consider the following distributed optimization over an undirected and connected network with $N$ agents.
\begin{align}\label{ProblemDIS}
\tilde{\m{x}}^*\in\mathop{\arg\min}\limits_{\tilde{\m{x}}\in \mathbb{R}^l}~\big\{\sum_{i=1}^N \big(f_i(\tilde{\m{x}})+r_i(\m{B}_i\tilde{\m{x}})\big)\big\},
\end{align}
where $f_i:\mathbb{R}^{l}\rightarrow(-\infty,+\infty)$ and $r_i:\mathbb{R}^{m_i}\rightarrow(-\infty,+\infty]$ are convex and accessed only by agent $i$, and $\m{B}_i\in \mathbb{R}^{m_i\times l}$ is a given matrix. Each agent $i$ aims to obtain an optimal solution $\tilde{\m{x}}^*$ of \eqref{ProblemDIS} by the local computation and communication with its immediate neighbors. Inspired by \cite{Guo2022}, this composite distributed optimization can be reformulated as
  $\min_{\m{x}}\sum_{i=1}^{N} (f_i(\m{x}_i)+r_i(\m{B}_i\m{x}_i))+\delta_0(\m{Dx}),$
  where $\m{x}_i\in\mathbb{R}^{l}$, $\m{x}=[(\m{x}_1)\tr,\cdots,(\m{x}_N)\tr]\tr$, and $\delta_0(\m{Dx})$ encodes the constraint $\m{x}_i=\m{x}_j, i,j=1,\cdots,N$. In this problem, $f_i(\m{x}_i)$ is a private loss function, and $r_i(\m{B_ix_i})$ is usually a regularizer. Therefore, it can be viewed as a special form of problem \eqref{Problem}.
%(general LASSO regularizer \cite{Example3}, fused LASSO regularizer \cite{Tibshirani2005}, group LASSO regularizer \cite{Jacob2009}, and OSCAR regularizer \cite{Howard2008}).
\begin{table*}[!t]
\renewcommand\arraystretch{1.5}
\begin{center}
\caption{The comparison of BALPA with existing work on PD-PSAs ($\alpha$ and $\beta$ are the primal and dual update stepsizes, respectively).}
\scalebox{1}{
\begin{tabular}{ccccc}
\hline
\textbf{Problem}&\textbf{Algorithm}&\textbf{Convergence Condition}&\textbf{\# of $\nabla f(\cdot)$ / Iter.}&\textbf{\# of $\mathrm{prox}^{\alpha}_{r}(\cdot)$ / Iter.}\\
\hline
$f=0$, $\m{B}=I$&L-ALM, PDHG&$0<\alpha\beta\|\m{D}\tr\m{D}\|<1$, $\alpha>0$&--&Once\\
$f\neq0$, $\m{B}=I$&L-ALM, C-V, TriPD&$0<\alpha\beta\|\m{D}\tr\m{D}\|+\alpha L/2<1$, $\alpha>0$&Once&Once\\
$f\neq0$, $\m{B}=I$&PDFP&$0<\alpha<2/L$, $0<\alpha\beta\|\m{D}\tr\m{D}\|<1$&Twice&Twice\\
$f\neq0$, $\m{B}=I$&PD3O&$0<\alpha<2/L$, $0<\alpha\beta\|\m{D}\tr\m{D}\|<1$&Twice&Once\\
$f\neq0$, $\m{B}=I$&AFBA, PDDY&$0<\alpha<2/L$, $0<\alpha\beta\|\m{D}\tr\m{D}\|<1$&Once&Once\\
\hline
\eqref{Problem}, TranA&C-V, TriPD&$0<\alpha\beta\|\m{B}\tr\m{B}+\m{D}\tr\m{D}\|+\alpha L/2<1$, $\alpha>0$&Once&Once\\
\eqref{Problem}, TranA&PAPC-type&$0<\alpha<2/L$, $0<\alpha\beta\|\m{B}\tr\m{B}+\m{D}\tr\m{D}\|<1$&Once&Once\\
\eqref{Problem}, TranB&L-ADMM&$0<\alpha\beta\|\m{B}\tr\m{B}+\m{D}\tr\m{D}\|+\alpha L/2<1$, $\alpha>0$&Once&Once\\
\eqref{Problem}, TranB&C-V, TriPD, L-LAM&$\alpha\beta\|\M{D}_1\M{D}_1\tr+\M{D}_2\M{D}_2\tr\|+\alpha L/2<1$, $\alpha>0$&Once&Once\\
\eqref{Problem}, TranB&PDFP&$0<\alpha<2/L$, $0<\alpha\beta\|\M{D}_1\M{D}_1\tr+\M{D}_2\M{D}_2\tr\|<1$&Twice&Twice\\
\eqref{Problem}, TranB&PD3O&$0<\alpha<2/L$, $0<\alpha\beta\|\M{D}_1\M{D}_1\tr+\M{D}_2\M{D}_2\tr\|<1$&Twice&Once\\
\eqref{Problem}, TranB&AFBA, PDDY&$0<\alpha<2/L$, $0<\alpha\beta\|\M{D}_1\M{D}_1\tr+\M{D}_2\M{D}_2\tr\|<1$&Once&Once\\
\hline
\eqref{Problem}&\textbf{BALPA}&$0<\alpha<2/L$&Once&Once\\
\hline
\end{tabular}}
\label{TableCOMPA}
\end{center}
\end{table*}

\subsection{Related Work}
To date, there have been numerous primal--dual proximal splitting algorithms (PD-PSAs) proposed for problem \eqref{Problem} and its special cases. Firstly, consider the problem $\min_{\m{x}}r(\m{x})+\delta_0(\m{Dx}-\m{d})$. To this special case, the linearized augmented Lagrangian method (L-ALM) \cite{LALM2,LALM3} and the primal--dual hybrid gradient algorithm (PDHG) \cite{LALM1,CP2} are provided under the convergence condition $0<\alpha\beta\|\m{D}\tr\m{D}\|<1$, where $\alpha>0$ and $\beta>0$ are the primal and dual update stepsizes, respectively. Then, consider the problem $\min_{\m{x}}f(\m{x})+r(\m{x})+\delta_0(\m{Dx}-\m{d})$. By linearizing the smooth term $f(\m{x})+\frac{\beta}{2}\|\m{Dx}-\m{d}\|^2$, \cite{ADMM-D1} proposes another L-ALM and establishes the $O(\frac{1}{k})$ convergence rate in the
ergodic sense when the stepsizes satisfy $\alpha\beta\|\m{D}\tr\m{D}\|+\alpha L<1$. In addition, the Condat-Vu algorithm (C-V) \cite{CV1,CV2} and the triangularly preconditioned primal--dual algorithm (TriPD) \cite{TriPD} are provided based on the forward-backward splitting (FBS) \cite{FBS} and asymmetric forward-backward splitting (AFBS), respectively, under the convergence condition $\alpha\beta\|\m{D}\tr\m{D}\|+{\alpha L}/{2}<1$. To allow for a larger range of acceptable stepsizes, some other PD-PSAs for minimizing the sum of three convex functions, such as the forward-backward-adjoint splitting algorithm (AFBA) \cite{AFBA}, the primal--dual fixed-point algorithm (PDFP) \cite{PDFP}, the primal--dual three-operator splitting (PD3O) \cite{PD3O}, and the Davis-Yin splitting \cite{DYS-T} based primal--dual algorithm (PDDY) \cite{DYS-T2}, are provided for this special case. The convergence conditions of these PD-PSAs are $0<\alpha<2/L$, $0<\alpha\beta\|\m{D}\tr\m{D}\|<1$.

For solving problem \eqref{Problem}, on one hand, we can reformulate the problem as $\min_{\m{x}} f(\m{x})+\tilde{h}([\m{B};\m{D}]\m{x})$, called TranA, where $\tilde{h}([\m{B};\m{D}]\m{x})=r(\m{B}\m{x})+\delta_0(\m{Dx}-\m{d})$. To this reformulation, some PD-PSAs for minimizing the sum of two convex functions can be applied, such as the Douglas-Rachford splitting dynamics \cite{DRSD}, the primal--dual fixed-point algorithm based on the proximity operator (PDFP$^2$O) \cite{PDFP2O}, and the proximal alternating predictor-corrector algorithm (PAPC) \cite{PAPC}. Besides, the PD-PSAs mentioned earlier, i.e., C-V, TriPD, AFBA, PDFP, PD3O, and PDDY, are also applicable to this reformulation. On the other hand, replacing $\m{Bx}$ by $\m{y}$, the problem \eqref{Problem} can be reformulated as $\min\{f(\m{x})+r(\m{y}):\M{D}_1\m{x}+\M{D}_2\m{y}=\M{d}\}$, where $\M{D}_1=[\m{D}\tr,\m{B}\tr]\tr$, $\M{D}_2=[0,-I]\tr$, and $\M{d}=[\m{d}\tr,0]\tr$, called TranB. To this reformulation, although the dimensionality of primal and dual variables is increased, both the loss function and the equation constraints are structured separable, which implies that it can be efficiently solved by L-ALM and the linearized alternating direction method of multipliers (L-ADMM) \cite{ADMM-D1,ADMM-D3}. Furthermore, the PD-PSAs proposed in \cite{CV1,CV2,AFBA,TriPD,PDFP,PD3O,DYS-T2} for minimizing the sum of three convex functions can also be used to solve this reformulation.

We summary the convergence conditions for these PD-PSAs under various situations in Table \ref{TableCOMPA}, where PDFP, AFBA, PD3O are called a class of PAPC-type algorithms, since they will degenerate to PAPC when minimizing $ f(\m{x})+\tilde{h}([\m{B};\m{D}]\m{x})$. Note that the convergence condition of all these algorithm depends on the linear mappings $\m{B}$ and $\m{D}$. It implies that when $\|\m{B}\tr\m{B}\|$ or $\|\m{D}\tr \m{D}\|$ is large, the stepsizes $\alpha$ and $\beta$ of corresponding algorithms are forced to be small to guarantee their convergence conditions. Therefore, more iterations are required, which may lead to a significant reduction in the efficiency of these algorithms. To overcome the dependence of convergence condition on linear mapping, the indefinite-proximal ALM \cite{InexactProximalALM1} and the indefinite-proximal ADMM \cite{InexactProximalALM1,InexactProximalALM2} are given. To these algorithms, they achieve a linear mapping-independent convergence condition at the cost of increasing the difficulty of primal updates. However, for the primal update of these algorithms, there is usually no explicit iterative scheme, but rather a subproblem to be solved. Therefore, for implementation, one has to design additional optimization algorithms to solve this series of subproblems (exactly or inexactly) based on the characteristics of the considered problem.

As an important application of problem \eqref{Problem}, let us turn our attention to distributed optimization problem \eqref{ProblemDIS}.
%To achieve the exact convergence under a fixed stepsize, many distributed primal--dual algorithms have been proposed in \cite{TriPD,Sulaiman2021,Xu2021,Chang2015,PGEXTR,NDIS,AMM2022,PDFP-DIS,TPUS}. In \cite{Sulaiman2021,Xu2021,Chang2015,PGEXTR,NDIS,AMM2022}, when $r_i\neq0$ and $\m{B}_i=I$, with the assumption that $r_i$ is proximable, distributed PD-PSAs have been presented. When $r_i\neq0$ and $\m{B}_i\neq I$, based on AFBS (or Tri-PD), a distributed Tri-PD (TriPD-Dist) has been provided in \cite{TriPD}.
Based on PDFP and PD3O, \cite{PDFP-DIS} and \cite{TPUS} provide the distributed PDFP (PDFP-Dist) and the triple proximal splitting algorithm with uncoordinated stepsizes (TPUS), respectively. Moreover, by inexact FBS and Fenchel-Moreau-Rockafellar duality, the dual inexact splitting algorithm (DISA) \cite{Guo2022} is provided, which is the first distributed optimization algorithm with $\|\m{B}_i\tr\m{B}_i\|$-independent convergence condition. However, two proximal mappings of $r_i$ are required in each iteration of DISA, which implies that the additional cost of one proximal mapping in DISA can not be ignored if $\mathrm{prox}^{\alpha}_{r_i}(\cdot)$ is difficult to compute.

\subsection{Our Contributions}
To overcome the dependence of the convergence condition on the linear mappings $\m{B}$ and $\m{D}$, we provide a novel PD-PSA called \underline{Bal}anced \underline{P}rimal--Dual Proximal \underline{A}lgorithm (BALPA). Without any additional condition to the considered problem \eqref{Problem}, the convergence condition of this new algorithm is
\begin{equation}\label{ConSTEP}
0<\alpha<2/L.
\end{equation}
This convergence condition is not restricted by any condition related to $\m{B}$ and $\m{D}$ explicitly or implicitly. Different from existing PD-PSAs, tiny stepsizes can be avoided even when $\|\m{B}\tr\m{B}\|$ and $\|\m{D}\tr\m{D}\|$ are large, thus ensuring a fast convergence of BALPA. Moreover, the stepsize selection of BALPA is more flexible, and there is no need to estimate $\|\m{B}\tr\m{B}\|$ and $\|\m{D}\tr\m{D}\|$ when tuning the stepsizes. We consider BALPA as a necessary supplement to DP-PSAs, especially for the case where $r$ is proximable but $\|\m{B}\tr\m{B}\|$ or $\|\m{D}\tr\m{D}\|$ is large.

To the development of BALPA, different from the indefinite-proximal method  \cite{InexactProximalALM1,InexactProximalALM2}, we balance the primal update and the dual update to eliminate the dependence of the convergence condition on the linear mapping. This idea comes mainly from the recently studied preconditioning technique for the PDHG \cite{Preconditioned2,BALMhe}. In each iteration of BALPA, to maintains the proximity-induced feature, the dual update is a proximal point of a specifically designed time-varying function $\psi_{k}$, whose proximal mapping has a closed-form representation, and the primal update is the same as the existing DP-PSAs but with an additional correction step. Therefore, when $\mathrm{prox}_{\psi_{k}}^{\gamma}(\cdot)$ is known, where $\gamma>0$ is an arbitrary constant, BALPA inherits the advantages of simple structure and easy implementation of the traditional DP-PSAs. Even though $\mathrm{prox}_{\psi_{k}}^{\gamma}(\cdot)$ is unknown, the dual update is equivalent to solving a positive definite system of linear equations, which can be easily solved by the preconditioning conjugate gradient method. Therefore, the cost of this new dual update is perfectly acceptable. This is exactly the distinctive merit of BALPA. Furthermore, we propose the stochastic BALPA (S-BALPA), where a random estimate of the gradient of $f$ is used, instead of the true gradient.

In addition, we apply the developed BALPA to distributed composite optimization problem \eqref{ProblemDIS}, thus obtaining a novel efficient distributed algorithm called BALPA-Dist. The convergence condition of the proposed algorithm is neither restricted by the network topology nor by $\|\m{B}_i\tr\m{B}\|$. We compare it with existing distributed PD-PSAs for problem \eqref{ProblemDIS} in Table \ref{TableCOMPA2}. As illustrated in Table \ref{TableCOMPA2}, BALPA-Dist has a more relaxed convergence condition. Moreover, in each iteration, $\nabla f_i(\cdot)$ and $\mathrm{prox}_{r_i}^{\alpha}(\cdot)$ are calculated only once.
\begin{table}[!t]
\renewcommand\arraystretch{1.5}
\begin{center}
\caption{The comparison of BALPA-Dist with existing distributed PD-PSAs for problem \eqref{ProblemDIS}.}
\scalebox{0.9}{
\begin{tabular}{ccccc}
\hline
\textbf{Algorithm}&\textbf{\# of $\nabla f_i(\cdot)$ / Iter.}&\textbf{\# of $\mathrm{prox}^{\alpha}_{r_i}(\cdot)$ / Iter.}&\textbf{LMI}&\textbf{NI}\\
\hline
TriDP-Dist&Once&Once&\XSolidBrush&\XSolidBrush\\
PDFP-Dist&Twice&Once&\XSolidBrush&\Checkmark\\
TPUS&Four times&Once&\XSolidBrush&\Checkmark\\
DISA&Once&Twice&\Checkmark&\Checkmark\\
\hline
\textbf{BALPA-Dist}&Once&Once&\Checkmark&\Checkmark\\
\hline
\multicolumn{5}{l}{{\textbf{LMI}: Linear mapping independent convergence condition.}}\\
\multicolumn{5}{l}{{\textbf{NI}: Network independent convergence condition.}}
\end{tabular}}
\label{TableCOMPA2}
\end{center}
\end{table}

For convergence analysis, with the convergence condition \eqref{ConSTEP}, we prove the convergence of BALPA and establish the $O(1/k)$ non-ergodic convergence rate. In addition, with a stronger convergence condition that $0<\alpha<{1}/{L}$, we establish the $O(1/k)$ ergodic convergence rate in the primal--dual gap, the optimality gap, and the constraint violation, respectively. Moreover, with the diminishing stepsizes which are not restricted by $\m{B}$ and $\m{D}$, the $O({1}/{\sqrt{k}})$ convergence rate of S-BALPA is established. Under the strong convexity of $f$, the rate can be improved to $O(\ln k/k)$. With the constant stepsizes (not restricted by $\m{B}$ and $\m{D}$) and variance-reduced method, we prove the $O(1/k)$ rate of S-BALPA.

\subsection{Notations and Organization}
Notations: $\mathbb{R}^n$ denotes the $n$-dimensional vector space with inner-product $\langle \cdot,\cdot\rangle$. The $1$-norm and Euclidean norm are denoted as $\|\cdot\|_1$ and $\|\cdot\|$, respectively. $0$ and $I$ denote the null matrix and the identity matrix of appropriate dimensions, respectively. Let $1_n\in \mathbb{R}^n$ be a vector with each component being one. For a matrix $A\in \mathbb{R}^{p\times n}$, $\|A\|:=\max_{\|x\|=1} \|Ax\|$. If $A$ is symmetric, $A\succ 0$ (resp. $A\succeq 0$) means that $A$ is positive definite (resp. positive semi-definite). For a given matrix $H \succeq0$, $H$-norm is defined as $\|x\|^2_{H}=\langle x,H x\rangle$.

Organization: In Section \ref{SEC:Algorithm Development}, we propose BALPA for the general composite optimization problem \eqref{Problem}. Then, we make a simple modification over BALPA with a non-diagonal preconditioner to obtain BALPA-Dist for the distributed optimization problem \eqref{ProblemDIS}. In addition, we compare the proposed algorithms with the existing algorithms. In Section \ref{SEC:Convergence}, we prove the convergence of the proposed algorithms under a unified algorithmic framework. Moreover, we establish the sublinear convergence rate of the proposed algorithm and their stochastic versions. In Sections \ref{SEC:Simulation} and \ref{SEC:Conclusion}, two numerical  experiments and conclusions are given, respectively.

\section{Algorithm Development}\label{SEC:Algorithm Development}
In this section, we propose a novel PD-PSA for problem \eqref{Problem}. Throughout this paper, we give the following assumption.
\begin{Assumption}\label{ASS1}
The set of solutions to problem \eqref{Problem} is nonempty. The function $f(\m{x})$ and $r(\m{x})$ are proper closed convex functions. Besides, $f(\m{x})$ is a smooth function with $L$-Lipschitz continuous gradient.
\end{Assumption}
With Assumption \ref{ASS1}, for any $\m{x}_1,\m{x}_2,\m{x}_3\in\mathbb{R}^{n}$, we have the following two commonly used inequalities
\begin{align}
&\langle \m{x}_1-\m{x}_2,\nabla f(\m{x}_3)-\nabla f(\m{x}_1)\rangle\leq \frac{L}{4} \|\m{x}_2-\m{x}_3\|^2,\label{SM1}\\
&\langle \m{x}_1-\m{x}_2,\nabla f(\m{x}_3)\rangle\leq f(\m{x}_1)-f(\m{x}_2)+\frac{L}{2}\|\m{x}_2-\m{x}_3\|^2.\label{SM2}
\end{align}

\subsection{Balanced Primal--Dual Algorithm}
Introducing the equality constraint $\m{Bx}=\m{y}$, and letting $\M{X}=[\m{x}\tr,\m{y}\tr]\tr\in\mathcal{X}:=\mathbb{R}^{n+p_1}$, $F(\M{X})=f(\m{x})+0,~R(\M{X})=0+r(\m{y})$, $\M{d}=[\m{d}\tr,0]\tr\in\mathcal{Y}:=\mathbb{R}^{p_2+p_1}$, and $\M{D}:\mathbb{R}^{n+p_1}\rightarrow \mathbb{R}^{p_2+p_1}:(\m{x},\m{y})\mapsto(\m{Dx},\m{Bx-y})$, problem \eqref{Problem} is equivalent to
\begin{equation}\label{ProblemNew}
\min_{\m{x}} (F+R)(\M{X})+\delta_0(\M{DX}-\M{d}),
\end{equation}
where $\delta_0(\M{DX}-\M{d})=0$ if $\M{DX}-\M{d}=0$; otherwise $\delta_0(\M{DX}-\M{d})=\infty$. By strongly duality, to solve problem \eqref{ProblemNew}, we can consider the saddle-point problem
\begin{align}\label{ProblemSaddle}
\min_{\M{X}}\max_{\M{\Lambda}}\mathcal{L}(\M{X},\M{\Lambda})=(F+R)(\M{X})+\langle \M{\Lambda},\M{DX-d} \rangle,
\end{align}
where $\M{\Lambda}\in \mathbb{R}^{p_2+p_1}$ is the Lagrange multiplier. Let $\alpha>0$ be a given constant, $\gamma>0$ be an arbitrary constant, $\{\alpha_k\}_{k\geq0}$ be the monotonically non-increasing sequence of the stepsize of primal update satisfying $0<\alpha_k\leq\alpha,k\geq0$, and $\M{g}^{k}=[(g^k)\tr,0]\tr$ be the stochastic gradient of $\nabla F(\M{X}^{k})$. To solve \eqref{ProblemSaddle}, we balance the primal and dual update and propose stochastic BALPA (S-BALPA):
\begin{flushleft}
\centering\fbox{
 \parbox{0.43\textwidth}{
\textbf{Algorithm 1: S-BALPA for problem \eqref{Problem}}.\\
Initialize $\M{X}^0\in \mathbb{R}^{n+p_1}$ and $\M{\Lambda}^0=0$. Define
$$\psi_{k}(\M{\Lambda})=\frac{\alpha}{2}\|\M{D}\tr(\bm{\Lambda}-\bm{\Lambda}^k)\|^2+\langle\bm{\Lambda},\M{d}-\M{D}\bar{\M{X}}^k\rangle.$$
With $(\M{X}^k,\M{\Lambda}^k)$, the new iterate $(\bar{\M{X}}^k,\M{X}^{k+1},\M{\Lambda}^{k+1})$ is generated via the
following steps:
\begin{subequations}
\begin{align}
\bar{\M{X}}^{k}&=\mathrm{prox}^{\alpha_{k}}_{R}(\M{X}^{k}-\alpha_{k}(\M{D}\tr\bm{\Lambda}^{k}+ \M{g}^{k})),\\
\bm{\Lambda}^{k+1}&=\mathrm{prox}^{\gamma}_{\psi_{k}}(\M{\Lambda}^k),\label{DUALBAS1}\\
\M{X}^{k+1}&=\bar{\M{X}}^{k}+\alpha_k\M{D}\tr(\bm{\Lambda}^k-\bm{\Lambda}^{k+1}).
\end{align}
\end{subequations}
}}\end{flushleft}
Different from the classic PD-PSAs, where the dual update is $\M{\Lambda}^{k+1}=\M{\Lambda}^k+\beta(\M{DX}^{k+1}-\M{d})$ or $\M{\Lambda}^{k+1}=\M{\Lambda}^k+\beta(\M{D}\bar{\M{X}}^{k}-\M{d})$, the implementation costs of the primal update and dual update of BALPA are balanced. More specifically, the primal update and dual update are the proximal mappings with closed-form representation. Therefore, we name the proposed algorithm as balanced primal--dual algorithm. The main features of BALPA (S-BALPA) are elaborated as follows.
\begin{itemize}
  \item \textbf{The recursion of BALPA is proximity-induced.} Both the primal update and the dual update of BALPA enjoy the proximity-induced feature.
  \item \textbf{The computation cost of BALPA is acceptable.} Since the time-varying function $\psi_k$ is quadratic and strongly convex, $\mathrm{prox}_{\psi_k}^{\gamma}(\cdot)$ has a unique closed-form representation. When the closed-form representation is known, the computational complexity is the same as the classic PD-PSAs. When $\mathrm{prox}_{\psi_k}^{\gamma}(\cdot)$ is unknown, by the definition of $\mathrm{prox}_{\psi_k}^{\gamma}(\cdot)$, it holds that
      $$\bm{\Lambda}^{k+1}=\mathop{\arg\min}\limits_{\bm{\Lambda}}\{\frac{1}{2}\|\bm{\Lambda}-\bm{\Lambda}^k\|^2_{\M{Q}}+\langle\bm{\Lambda},\M{d}-\M{D}\bar{\M{X}}^k\rangle\},$$
where $\mathbf{Q}=\frac{1}{\gamma}I+\alpha\M{DD}\tr$. Hence, it can be solved effectively by the preconditioning conjugate gradient method.
  \item \textbf{The convergence condition of BALPA is independent of the linear mappings.} By balancing the primal and dual update, the profit mapping of BALPA is
  $$
\M{M}:(\M{X},\M{\Lambda})\mapsto ((\frac{1}{\alpha}-\frac{L}{2})\M{X},\frac{1}{\gamma}\M{\Lambda}).
  $$
  In the next section, we will show that the positive definiteness of the profit mapping $\M{M}$ ensures the convergence of BALPA. Therefore, the convergence condition of BALPA is $0<\alpha<2/L$, which is not restricted by any condition related to B and D explicitly or implicitly.
\end{itemize}
To summarize, BALPA not only retains the advantages of simple structure and easy implementation of the classic PD-PSAs, but also overcomes the high dependence of stepsize selection of PD-PSAs on linear mapping.
\subsection{Compared to Existing PD-PSAs}
This subsection compares BALPA to existing PD-PSAs, such as C-V, TriPD, PD3O, PDFP, and AFBS.

We start to compare the convergence conditions of these existing PD-PSAs with BALPA. Recall profit mappings of these existing PD-PSAs. The profit mappings of C-V and TriPD are
  $$
\M{M}_{\m{CV}}:(\M{X},\M{\Lambda})\mapsto ((\frac{1}{\alpha}-\frac{L}{2})\M{X}-\M{D}\tr\M{\Lambda},-\M{DX}+\frac{1}{\beta}\M{\Lambda}),
  $$
and the profit mappings of PD3O, PDFP, and AFBS are
  $$
\M{M}_{\m{LV}}:(\M{X},\M{\Lambda})\mapsto ((\frac{1}{\alpha}-\frac{L}{2})\M{X},(\frac{1}{\beta}-\alpha\M{DD}\tr)\M{\Lambda}).
  $$
To ensure the positive definiteness of $\M{M}_{\m{CV}}$ and $\M{M}_{\m{LV}}$, the convergence conditions of C-V and TriPD are
$$\alpha>0,\beta>0,\alpha\beta\|\M{D}\tr\M{D}\|+\alpha L/2<1,$$ and the stepsize conditions of PD3O, PDFP, and AFBA are
$$0<\alpha<2/L, \alpha\beta\|\M{D}\tr\M{D}\|<1.$$
In BALPA, the convergence condition is $0<\alpha<2/L$, which is independent of the mappings $\m{B}$ and $\m{D}$.

Then, we compare the iteration of BALPA with the existing PD-PSAs for problem \eqref{ProblemNew} in Table \ref{Table-Iteration}, where using $(\M{X},\M{\Lambda})$ and $(\M{X}^+,\M{\Lambda}^+)$ for the current and next iterations. As shown in Table \ref{Table-Iteration}, the primal updates of these PD-PSAs are identical, with the main difference being the dual updates and the corrections. More specifically, C-V, TriPD, and PD3O are correction-then-dual-update algorithms. PDFP, AFBS, and BALPA are dual-update-then-correction algorithms. In the correction of PD3O, to achieve a more relaxed convergence condition, there are two more terms $(\nabla F(\M{X})$ and $\nabla F(\bar{\M{X}}))$ than C-V and TriPD in correction. In fact, $\nabla F(\M{X})$ has been computed in the previous step, and only $\nabla F(\bar{\M{X}})$ is required to be computed in current iteration. Therefore, there is no additional cost in PD3O compared to C-V and TriPD, except that $\nabla F(\M{X})$ needs to be stored. In the correction of PDFP, $\mathrm{prox}_{R}^{\alpha}(\cdot)$ is applied, and it will not be used in the next iterations. In the correction of AFBA and BALPA, only the mapping $\M{D}\tr$ is needed, thus it has the advantage of low computational and storage costs.

\begin{table}[!t]
\renewcommand\arraystretch{1.5}
\begin{center}
\caption{The iterations of BALPA and existing PD-PSAs for problem \eqref{ProblemNew}.}
\scalebox{1}{
\begin{tabular}{|c|l|}
\hline
\multicolumn{2}{|c|}{$\bar{\M{X}}=\mathrm{prox}^{\alpha}_{R}(\M{X}-\alpha(\M{D}\tr\bm{\Lambda}+ \nabla F(\M{X}))$}\\
\hline
\textbf{Algorithm}&\textbf{Dual Update \& Correction}\\
\hline
\hline
\multirow{2}{*}{C-V, TriPD}     & $\M{X}^+=2\bar{\M{X}}-\M{X}$\\
                                         & $\bm{\Lambda}^+=\M{\Lambda}+\beta(\M{D}\M{X}^+-\M{d})$\\
\hline
\multirow{2}{*}{PD3O}     & $\M{X}^+=2\bar{\M{X}}-\M{X}+\alpha(\nabla F(\M{X})-\nabla F(\bar{\M{X}}))$\\
                                         & $\M{\Lambda}^+=\M{\Lambda}+\beta(\M{D}\M{X}^+-\M{d})$\\
\hline
\multirow{2}{*}{PDFP}     & $\bm{\Lambda}^+=\M{\Lambda}+\beta(\M{D}\bar{\M{X}}-\M{d})$\\
                                         & $\M{X}^+=\mathrm{prox}^{\alpha}_{R}(\M{X}-\alpha(\M{D}\tr\bm{\Lambda}^++ \nabla F(\M{X}))$\\
\hline
\multirow{2}{*}{AFBA}     & $\bm{\Lambda}^+=\M{\Lambda}+\beta(\M{D}\bar{\M{X}}-\M{d})$\\
                                         & $\M{X}^+=\bar{\M{X}}+\alpha\M{D}\tr(\bm{\Lambda}-\bm{\Lambda}^+)$\\
\hline
\multirow{2}{*}{\textbf{BALPA}}     & $\bm{\Lambda}^+=\arg\min_{\M{s}}\big\{\frac{1}{2}\|\M{s}-\bm{\Lambda}\|^2_{\M{Q}}+\langle\M{s},\M{d}-\M{D}\bar{\M{X}}\rangle\big\}$\\
                                         & $\M{X}^+=\bar{\M{X}}+\alpha\M{D}\tr(\bm{\Lambda}-\bm{\Lambda}^+)$\\
\hline
\end{tabular}}
\label{Table-Iteration}
\end{center}
\end{table}

\subsection{BALPA for Distributed Optimization Problems}
Consider the distributed optimization problem \eqref{ProblemDIS} on an undirected and connected network $\mathcal{G}(\mathcal{V},\mathcal{E})$. Similar as \cite{Guo2022}, we introduce the local copy $\m{x}_i\in \mathbb{R}^l$ of $\tilde{\m{x}}$ and a symmetric and doubly stochastic matrix $U\in\mathbb{R}^{N\times N}$ satisfying that $\mathrm{null}(I-U)=\mathrm{span}(1_N)$, and if $(i,j)\notin \mathcal{E}$ and $i\neq j$, $U_{ij}=U_{ji}=0$, otherwise $U_{ij}>0$. Let $\m{x}=[(\m{x}_1)\tr,\cdots,(\m{x}_N)\tr]\tr$, $\m{B}$ be a block diagonal matrix with its $(i,i)$-th block being $\m{B}_i$ and other blocks being $0$, and $\m{D}=(\frac{1}{2}(I-U)\otimes I_l)^{\frac{1}{2}}$. The distributed optimization problem \eqref{ProblemDIS} is equivalent to
$\min_{\m{x}} f(\m{x})+r(\m{B}x)+\delta_0(\m{Dx})$, where $f(\m{x})=\sum_{i=1}^{N}f_i(\m{x}_i)$, and $r(\m{Bx})=\sum_{i=1}^{N}r_i(\m{B}_i\m{x}_i)$. Therefore, it can be solved by BALPA (S-BALPA). However, the dual update \eqref{DUALBAS1} of BALPA (S-BALPA) for problem \eqref{ProblemDIS} cannot be implemented in a distributed manner. We need to further improve it technically.

For problem \eqref{ProblemDIS}, the dual update of BALPA is equivalent to
$
\bm{\Lambda}^{k+1}=\arg\min_{\bm{\Lambda}}\{\frac{1}{2}\|\bm{\Lambda}-\bm{\Lambda}^k\|^2_{\M{Q}}-\langle\bm{\Lambda},\M{D}\bar{\M{X}}^k\rangle\},
$
where $\mathbf{Q}=\frac{1}{\gamma}I+\alpha\M{DD}\tr$ and can be seen as a preconditioner. Therefore, inspired by \cite{Guo2022}, we can redesign the preconditioner for the dual update as follows.
$$
\M{Q}=
\left(
  \begin{array}{cc}
    \frac{\alpha}{\gamma}I & 0 \\
    0 & \m{S} \\
  \end{array}
\right),\text{ where } \m{S}=\frac{\alpha+\alpha\gamma}{\gamma}I+\frac{\alpha}{1-\gamma} \m{BB}\tr.
$$
By this non-diagonal preconditioning technique, we can directly derive the distributed BALPA (BALPA-Dist) from Algorithm 1. We elaborate the implementation
of Stochastic BALPA-Dist (S-BALPA-Dist) in Algorithm 2, where $\m{g}^k_i$ is the stochastic gradient of $\nabla f_i(\m{x}_i^k)$, and $\mu,\nu$ are the dual variables, which are associated with the equality constraints $\m{D}\m{x}=0$ and $\m{Bx}=\m{y}$, respectively. In the next section (Section \ref{SEC:Convergence}), we will show that the positive definiteness of $\M{Q}-\alpha\M{DD}\tr$ is an important prerequisite for the convergence of BALPA-Dist. According to \cite{Guo2022}, we have $\M{Q}-\alpha\M{DD}\tr\succ0$ when $\gamma\in(0,1)$, which implies that the positive definition of $\M{Q}-\alpha\M{DD}\tr$ is not restricted by $\|\m{B}_i\tr\m{B}_i\|$. Therefore, resorting to this non-diagonal preconditioner, BALPA-Dist can be implemented in a distributed manner, with a linear-mapping independent stepsize selection. For communication cost of BALPA-Dist, there is only one round of communication in each iteration. For computation costs, BALPA-Dist inherits the advantages of simple structure and easy implementation of BALPA.

\begin{flushleft}
\centering\fbox{
 \parbox{0.43\textwidth}{
\textbf{Algorithm 2: S-BALPA-Dist for \eqref{ProblemDIS}}.\\
Initialize $(\m{x}^0,\m{y}^0)\in \mathbb{R}^{n+p_1}$ and $(\mu^0,\nu^0)=0$. Define
$$\psi^i_{k}(\nu_i)=\frac{\alpha}{2(1-\gamma)}\|\m{B}_i\tr(\nu_i-\nu_i^k)\|^2-\langle\nu_i,\m{B}_i\bar{\m{x}}^k-\bar{\m{y}}^k_i\rangle.$$
With $(\m{x}^k,\m{y}^k,\mu^k,\nu^k)$, the new iteration is generated via the
following steps:
\begin{align*}
&\bar{\m{x}}_i^k=\m{x}_i^k-\alpha_k(\mu_i^k+\m{B}_i\tr\nu_i^k+\m{g}^k_i),\\
&\bar{\m{y}}_i^k=\mathrm{prox}_{r_i}^{\alpha_k}(\m{y}_i^k+\alpha_k\nu_i^k),\\
&\mu_i^{k+1}=\mu_i^{k}+\frac{\gamma}{2\alpha}(\bar{\m{x}}_i^k- \sum\nolimits_{j=1}^{N}U_{ij}\bar{\m{x}}_j^k),\\
&\nu_i^k=\mathrm{prox}_{\psi^i_{k}}^{\gamma/(\alpha+\alpha\gamma)}(\nu_i^k),\\
&\m{x}_i^{k+1}=\bar{\m{x}}_i^k+\alpha_k(\mu_i^{k}-\mu_i^{k+1}+\m{B}_i\tr(\nu_i^k-\nu_i^{k+1})),\\
&\m{y}_i^{k+1}=\bar{\m{y}}_i^k-\alpha_k(\nu_i^k-\nu_i^{k+1}).
\end{align*}
}}\end{flushleft}
Next, we compare BALPA-Dist with DISA. In terms of iterative format, both BALPA-Dist and DISA are ATC (Adapt-Then-Combine)-based distributed algorithms and their recursions are proximity-induced. Besides, the convergence conditions of BALPA-Dist and DISA are $0<\alpha<\frac{2}{L}$, $0<\gamma<1$, which is independent of the network
topology and $\|\m{B}_i\tr\m{B}\|$. The main difference between BALPA-Dist and DISA is the correction step. Although DISA is not derived from primal--dual update with a correction, it can be interpreted as the primal--dual algorithm with a correction by equivalence transformation. More specifically, the correction step of DISA and BALPA-Dist, respectively given by
\begin{align*}
[\mathrm{DISA}]:~&\M{X}^+=\mathrm{prox}^{\alpha}_{R}(\M{X}-\alpha(\M{D}\tr\bm{\Lambda}^++ \nabla F(\M{X})),\\
[\mathrm{BALPA}]:~&\M{X}^+=\bar{\M{X}}+\alpha\M{D}\tr(\bm{\Lambda}-\bm{\Lambda}^+).
\end{align*}
It can be seen that DISA uses the PDFP-type correction. In DISA, two proximal mappings of $r_i$ are needed in each iteration. If the mapping is not easy to calculate, the additional cost of one proximal mapping cannot be ignored. However, in BALPA-Dist, only the mapping $\m{B}_i\tr$ is required. Therefore, BALPA-Dist not only retains the advantages of DISA, but also has lower computational costs than DISA.

\section{Convergence Analysis}\label{SEC:Convergence}
In this section, we establish the convergence and convergence rate of BALPA and BALPA-Dist and their stochastic versions under a unified framework.

It is clear from the previous analysis that the compact form of BALPA-Dist is the same as BALPA, and they can both be written in the following form
\begin{align}
\bar{\M{X}}^{k}&=\mathrm{prox}^{\alpha_{k}}_{R}(\M{X}^{k}-\alpha_{k}(\M{D}\tr\bm{\Lambda}^{k}+ \M{g}^{k})), \label{BAS3}\\
\bm{\Lambda}^{k+1}&=\mathop{\arg\min}\limits_{\bm{\Lambda}}\{\frac{1}{2}\|\bm{\Lambda}-\bm{\Lambda}^k\|^2_{\M{Q}}+\langle\bm{\Lambda},\M{d}-\M{D}\bar{\M{X}}^k\rangle\},\label{BAS1}\\
\M{X}^{k+1}&=\bar{\M{X}}^{k}+\alpha_k\M{D}\tr(\bm{\Lambda}^k-\bm{\Lambda}^{k+1}), \label{BAS2}
\end{align}
where $\{\alpha_k\}_{k\geq0}$ is a monotonically non-increasing stepsize sequence of primal update satisfying $0<\alpha_k\leq\alpha,k\geq0$. Therefore, we only need to analyze the convergence of this unified framework. Denote the set of saddle-points of \eqref{ProblemSaddle} as $\mathcal{W}^*$. It holds that $\mathcal{W}^*=\mathcal{X}^*\times\mathcal{Y}^*$, where $\mathcal{X}^*$ is the optimal solution set to problem \eqref{ProblemNew} and $\mathcal{Y}^*$ is the optimal solution set to its dual problem. Consider the saddle subdifferential
\begin{align*}
\partial \mathcal{L}=\underbrace{\left(
                               \begin{array}{c}
                                 \partial R(\M{X}) \\
                                 0 \\
                               \end{array}
                             \right)}_{:=\bm{A}(\M{X},\bm{\Lambda})}+\underbrace{\left(
                       \begin{array}{c}
                         \mathbf{D}\tr\bm{\Lambda} \\
                         -\M{DX}+\mathbf{d} \\
                       \end{array}
                     \right)}_{:=\bm{B}(\M{X},\bm{\Lambda})}+\underbrace{\left(
                                       \begin{array}{c}
                                         \nabla F(\M{X}) \\
                                         0 \\
                                       \end{array}
                                     \right)}_{:=\bm{C}(\M{X},\bm{\Lambda})}.
\end{align*}
By the first-order optimality condition, one has that $(\M{X}^*,\M{\Lambda}^*)\in\mathcal{W}^*$ if and only if $0\in\partial\mathcal{L}(\M{X}^*,\M{\Lambda}^*)$. Therefore, it is not difficult to check that when $\alpha_k\equiv\alpha$ and $\M{g}^k=\nabla F(\M{X}^k)$, any fixed point $(\M{X}^*,\M{\Lambda}^*,\bar{\M{X}}^*)$ of BALPA satisfies $\M{X}^*=\bar{\M{X}}^*$ and $(\M{X}^*,\M{\Lambda}^*)\in\mathcal{W}^*$. Denote $\mathcal{W}=\mathcal{X}\times \mathcal{Y}$ and $\M{W}=[\M{X}\tr,\M{\Lambda}\tr]\tr\in\mathcal{W}$. A point $\M{W}^*=[(\M{X}^*)\tr,(\M{\Lambda}^*)\tr]\tr$ is called a saddle point of problem \eqref{ProblemSaddle}, if
$
\mathcal{L}(\M{X}^{*}, \M{\Lambda}) \leq \mathcal{L}(\M{X}^{*}, \M{\Lambda}^{*}) \leq \mathcal{L}(\M{X}, \M{\Lambda}^{*}), \forall \M{W}\in \mathcal{W},
$
which can be alternatively rewritten as, for $\forall \mathbf{W}\in \mathcal{W}$,
\begin{align}\label{IV2}
R(\M{X})-R(\M{X}^*)+\langle\mathbf{W}-\mathbf{W}^*,(\bm{B}+\bm{C})\mathbf{W}^*\rangle\geq 0.
\end{align}
Let $\Phi(\M{X})=F(\M{X})+R(\M{X})$. Formula \eqref{IV2} is equivalent to
$
\mathcal{L}(\M{X},\M{\Lambda}^*)-\mathcal{L}(\M{X}^*,\M{\Lambda})=\Phi(\M{X})-\Phi(\M{X}^*)+\langle\mathbf{W}-\mathbf{W}^*,\bm{B}\mathbf{W}^*\rangle\geq 0$.  Denote the set of saddle-points of \eqref{ProblemSaddle} as $\mathcal{W}^*$. It holds that $\mathcal{W}^*=\mathcal{X}^*\times\mathcal{Y}^*$, where $\mathcal{X}^*$ is the optimal solution set to problem \eqref{ProblemNew} and $\mathcal{Y}^*$ is the optimal solution set to its dual problem. To simplify the notations in the later analysis, we define $\bm{\delta}^k\triangleq \M{g}^k-\nabla F(\M{X}^k)$, and define the following self-adjoint linear operators, which are important in the convergence analysis. Let $\M{H}_k: \mathcal{W}\rightarrow \mathcal{W}:(\M{X},\M{\Lambda})\mapsto ({1}/{\alpha_k}\M{X},\M{Q}\M{\Lambda})$,
\begin{align*}
\M{M}_k: \mathcal{W}\rightarrow \mathcal{W}:&(\M{X},\M{\Lambda})\mapsto ((\tau_k-\frac{L}{2})\M{X},(\M{Q}-\alpha_k\M{DD}\tr)\M{\Lambda}),\\
\M{G}_k: \mathcal{W}\rightarrow \mathcal{W}:&(\M{X},\M{\Lambda})\mapsto ((\tau_k-L)\M{X},(\M{Q}-\alpha_k\M{DD}\tr)\M{\Lambda}),
\end{align*}
where $\tau_k=\frac{1}{\alpha_k}-\mathbf{1}(\bm{\delta}^k)\eta_k$ with $\mathbf{1}(\bm{\delta}^k)=1$, if $\|\bm{\delta}^k\|>0$; $\mathbf{1}(\bm{\delta}^k)=0$, otherwise. We give the following standard assumption on the stochastic gradient.
\begin{Assumption}\label{Stochastic gradient}
The stochastic gradient $\M{g}^{k}$ of $\nabla F(\M{X}^k)$ is unbiased, i.e., $\mathbb{E}\M{g}^k=\nabla F(\M{X}^{k})$, and the stochastic gradient variance is bounded, i.e., there exists $\sigma>0$ such that $\mathbb{E}\|\M{g}^k-\nabla F(\M{X}^{k})\|^2\leq\sigma^2$.
\end{Assumption}
Then, to establish the global convergence, the following lemma is provided.
\begin{Lemma}\label{LEMADA1}
Suppose that Assumption \ref{ASS1} holds. Let $\M{V}^{k+1}=[(\bar{\M{X}}^{k})\tr,(\M{\Lambda}^{k+1})\tr]\tr$ and $\{\M{\Lambda}^{k},\M{X}^{k},\bar{\M{X}}^k\}_{k\geq0}$ be the iterate sequence generated by \eqref{BAS3}--\eqref{BAS2}. If $\{\alpha_k\}_{k\geq0}$ is chosen such that $\M{M}_k\succ0$, then
\begin{align}\label{FunA1}
&2\big(R(\bar{\M{X}}^k)-R(\M{X})+\big\langle\mathbf{V}^{k+1}-\mathbf{W},(\bm{B}+\bm{C})\mathbf{W}\big\rangle\big)\nonumber\\
&\leq 2\langle \M{X}-\M{X}^k,\bm{\delta}^k\rangle+\frac{1}{\eta_k}\|\bm{\delta}^k\|^2-\|\M{V}^{k+1}-\M{W}^k\|^2_{\M{M}_k}\nonumber\\
&+\|\M{W}-\M{W}^k\|^2_{\M{H}_k}-\|\M{W}-\M{W}^{k+1}\|^2_{\M{H}_k}, \forall \M{W}\in \mathcal{W},
\end{align}
where $\eta_k>0$ is any prescribed sequence. Moreover, if $\{\alpha_k\}_{k\geq0}$ is chosen such that $\M{G}_k\succ0$, we have
\begin{align}\label{FunA2}
&2\big(\Phi(\bar{\M{X}}^k)-\Phi(\M{X})+\big\langle\mathbf{V}^{k+1}-\mathbf{W},\bm{B}\mathbf{W}\big\rangle\big)\nonumber\\
&\leq 2\langle \M{X}-\M{X}^k,\bm{\delta}^k\rangle+\frac{1}{\eta_k}\|\bm{\delta}^k\|^2-\|\M{V}^{k+1}-\M{W}^k\|^2_{\M{G}_k}\nonumber\\
&+\|\M{W}-\M{W}^k\|^2_{\M{H}_k}-\|\M{W}-\M{W}^{k+1}\|^2_{\M{H}_k}, \forall \M{W}\in \mathcal{W}.
\end{align}
\end{Lemma}
{\it Proof}:
See Appendix A.
\hfill{$\blacksquare$}
\subsection{Convergence Analysis of BALPA and BALPA-Dist}
Let $\M{J}=\M{Q}-\alpha\M{DD}\tr$. We give the following theorem to show that the sequence $\{\M{X}^k,\M{\Lambda}^k\}$ generated by the unified framework \eqref{BAS3}--\eqref{BAS2} converges to a primal--dual optimal solution of problem \eqref{ProblemSaddle}, when $\sigma=0$.
\begin{Theorem}\label{TH1}
Suppose that Assumption \ref{ASS1} holds, and $\sigma=0$, i.e., $\M{g}^k=\nabla F(\M{X}^k)$ with probability 1 for any $\M{X}$. If $\alpha_k\equiv\alpha\in(0,2/L)$ and $\M{J}\succ0$, it holds that for any $\M{W}^*\in\mathcal{W}^*$, the sequence $\{\|\M{W}^k-\M{W}^*\|\}_{k\geq0}$ is bounded. Moreover, there exists $\M{W}^{\infty}\in \mathcal{W}^*$ such that $\lim_{k\rightarrow \infty}\M{W}^k=\M{W}^{\infty}$.
\end{Theorem}
{\it Proof}:
See Appendix B.
\hfill{$\blacksquare$}

For BALPA where $\mathbf{Q}=\frac{1}{\gamma}I+\alpha\M{DD}\tr$, we know that $\M{J}=\M{Q}-\alpha\M{DD}\tr=\frac{1}{\gamma}I\succ0$ if $\gamma>0$. Therefore, to ensure the convergence of BALPA, the convergence condition is $0<\alpha<2/L$, $\gamma>0$. For BALPA-Dist, it follows from Section 2.3 that $\M{J}\succ 0$ if $0<\gamma<1$. Thus, the convergence condition of BALPA-Dist is $0<\alpha<2/L$, $0<\gamma<1$.

Then, with general convexity condition, we prove the $o(1/k)$ non-ergodic convergence rate and the $O(1/k)$ ergodic convergence rate of the unified framework \eqref{BAS3}--\eqref{BAS2}, respectively.
\begin{Theorem}\label{TH:RATE1}
Under the same settings of Theorem \ref{TH1}, if $\alpha_k\equiv\alpha\in(0,2/L)$ and $\M{J}\succ0$, for any $\M{W}^*\in\mathcal{W}^*$
  \begin{subequations}\label{D-BPDA-RATA1}
  \begin{align}
  \lim_{k\rightarrow\infty}k\|\M{V}^{k+1}-\M{W}^k\|^2_{\M{M}}&=0,\\
  \min_{0\leq k\leq K}\{\|\M{V}^{k+1}-\M{W}^k\|^2_{\M{M}}\}&\leq\frac{\|\M{W}^0-\M{W}^*\|^2_{\M{H}}}{K+1},
  \end{align}
  \end{subequations}
where $\M{H}: \mathcal{W}\rightarrow \mathcal{W}:(\M{X},\M{\Lambda})\mapsto (\frac{1}{\alpha}\M{X},\M{Q}\M{\Lambda})$ is a self-adjoint linear mapping. Moreover, let $\bar{\M{X}}_K=\frac{1}{K}\sum_{k=0}^{K-1}\bar{\M{X}}^k$ and ${\M{\Lambda}}_K=\frac{1}{K}\sum_{k=0}^{K-1}{\M{\Lambda}}^{k+1}$.
If $\alpha_k\equiv\alpha\in(0,1/L)$ and $\M{J}\succ 0$, it holds that
\begin{align}\label{PDGAP}
\mathcal{L}(\bar{\M{X}}_K,\M{\Lambda})-\mathcal{L}(\M{X},{\M{\Lambda}}_K)\leq\frac{\|\M{W}^0-\M{W}\|^2_{\M{H}}}{2K},
\end{align}
for $\forall \M{W}\in \mathcal{W}$. Furthermore, letting $\rho=\max\{1+\|\M{\Lambda}^*\|,2\|\M{\Lambda}^*\|\}$ with $\M{\Lambda}^*\in \mathcal{Y}^*$, it holds that
\begin{subequations}\label{D-BPDA-RATA2}
\begin{align}
|\Phi(\bar{\M{X}}_K)-\Phi(\M{X}^*)|&\leq\frac{\frac{1}{\alpha}\|\M{X}^0-\M{X}^*\|^2+\rho^2\|\M{Q}\|}{2K},\label{D-BPDA-RATA211}\\
\|\M{D}\bar{\M{X}}_K-\M{d}\|&\leq\frac{\frac{1}{\alpha}\|\M{X}^0-\M{X}^*\|^2+\rho^2\|\M{Q}\|}{2K}.\label{D-BPDA-RATA222}
\end{align}
\end{subequations}
\end{Theorem}
{\it Proof}:
See Appendix C.
\hfill{$\blacksquare$}

It follows from \eqref{D-BPDA-RATA1} and $\sum_{k=0}^{\infty}\|\M{V}^{k+1}-\M{W}^k\|^2_{\M{M}}<\infty$ that the rate of the running-average successive difference of BALPA and BALPA-Dist is of order $O(1/k)$, and the running-best successive difference $\min_{0\leq k\leq K}\{\|\M{V}^{k+1}-\M{W}^k\|^2_{\M{M}}\}$ is of order $o(1/k)$. With a stronger convergence condition that $0<\alpha<1/L$, the primal--dual gap \eqref{PDGAP} is established. Furthermore, the ergodic $O(1/k)$ convergence rate of the optimality gap and the constraint violation are presented in \eqref{D-BPDA-RATA211} and \eqref{D-BPDA-RATA222}, respectively.

\subsection{Convergence Analysis of S-BALPA and S-BALPA-Dist}
We now present the $O(1/\sqrt{k})$ convergence rate of the unified framework \eqref{BAS3}--\eqref{BAS2} with diminishing stepsize.
\begin{Theorem}\label{TH2}
Suppose that Assumptions \ref{ASS1} and \ref{Stochastic gradient} hold. Let $(\M{X}^*,\M{\Lambda}^*)\in \mathcal{W}$, $\mathcal{D}(\M{X})=\sup_{\M{X}^*\in\mathcal{X}^*}\|\M{X}-\M{X}^*\|$ and $\mathcal{D}=\max\{\mathcal{D}(\M{X}^k):k=0,1,\cdots\}$. For any $K\geq 1$, choosing
\begin{align*}
\left\{\begin{array}{ll}
\eta_k=1+\sqrt{k},\alpha_k=\frac{1}{c+\sqrt{k}},& \mathcal{D}<\infty; \\
\eta_k=1+\sqrt{K-1},\alpha_k=\frac{1}{c+\sqrt{K-1}},& \mathcal{D}=\infty.
\end{array}
\right.
\end{align*}
If $c>1+\frac{L}{2}$ and $\M{J}\succ0$, which implies that $\M{M}_k\succ0$, one has
\begin{align}\label{SRATE1}
&\min_{0\leq k\leq K-1}\big\{\mathbb{E}\big(\|\M{V}^{k+1}-\M{W}^k\|^2_{\M{M}_k}\big)\big\}\leq\frac{\mathcal{C}_1}{\sqrt{K}}+\frac{\mathcal{C}_2}{K},
\end{align}
where $\mathcal{C}_1=\bar{\mathcal{D}}^2+\sigma^2$ with $\bar{\mathcal{D}}=\mathcal{D}$ when $\mathcal{D}<\infty$ and $\bar{\mathcal{D}}=\|\M{X}^0-\M{X}^*\|$ when $\mathcal{D}=\infty$, and $\mathcal{C}_2=c\|\M{X}^0-\M{X}^{*}\|^2+\|\M{\Lambda}^*\|^2_{\M{Q}}$. Moreover, if $c>1+L$ and $\M{J}\succ0$, which implies that $\M{G}_k\succ0$, one has that
\begin{align}\label{SRATE2}
2\mathbb{E}\big(\Phi(\bar{\M{X}}_K)-\Phi(\M{X}^*)+\rho\|\M{D}\bar{\M{X}}_K-\M{d}\|\big)\leq \frac{\mathcal{C}_1}{\sqrt{K}}+\frac{\mathcal{C}_3}{K},
\end{align}
where $\mathcal{C}_3=c\|\M{X}^0-\M{X}^{*}\|^2+\rho^2\|\M{Q}\|$.
\end{Theorem}
{\it Proof}:
See Appendix D.
\hfill{$\blacksquare$}

It follows from \eqref{SRATE1} that the running-best successive difference $\min_{0\leq k\leq K}\{\|\M{V}^{k+1}-\M{W}^k\|^2_{\M{H}_k}\}$ is of order $O(1/\sqrt{k})$. Moreover, similar as Theorem \ref{TH:RATE1}, assuming the existence of the dual optimal solution $\M{\Lambda}^*$, we can further assess the feasibility violation of the possibly infeasible solution $\bar{\M{X}}_K$ as in \eqref{SRATE2}. By Jensen's inequality it follows immediately that
$2(\Phi(\mathbb{E}(\bar{\M{X}}_K))-\Phi(\M{X}^*)+\rho\|\M{D}\mathbb{E}(\bar{\M{X}}_K)-\M{d}\|)\leq \frac{\mathcal{C}_1}{\sqrt{K}}+\frac{\mathcal{C}_3}{K}.$
Then, similar as the proof of Theorem \ref{TH:RATE1}, setting $\rho=\max\{1+\|\M{\Lambda}^*\|,2\|\M{\Lambda}^*\|\}$ and applying \cite[Lemma 2.3]{ADMM-D1} in this inequality, it holds that
\begin{subequations}\label{SRATE3}
\begin{align}
2|\Phi(\mathbb{E}(\bar{\M{X}}_K))-\Phi(\M{X}^*)|&\leq \frac{\mathcal{C}_1}{\sqrt{K}}+\frac{\mathcal{C}_3}{K},\\
2\|\M{D}(\mathbb{E}(\bar{\M{X}}_K))-\M{d}\|&\leq \frac{\mathcal{C}_1}{\sqrt{K}}+\frac{\mathcal{C}_3}{K}.
\end{align}
\end{subequations}
The same logic applies to all the subsequent convergence rate results.

With strong convexity of $f$, the convergence rate for the unified framework \eqref{BAS3}--\eqref{BAS2} with diminishing stepsize can be improved to $O(\ln k/k)$. The result is presented as follows.
\begin{Corollary}\label{COR:RATE1}
Suppose that Assumptions \ref{ASS1} and \ref{Stochastic gradient} hold, and $f$ is $\mu$-strongly convex, i.e., $f(\m{y})\geq f(\m{x})+\langle \nabla f(\m{x}),\m{y-x} \rangle+\frac{\mu}{2}\|\m{x-y}\|^2, \forall \m{x,y}$ with $\mu>0$. Choose $\eta_k=\mu (k+1)$ and $\frac{1}{\alpha_k}=c+\mu(k+1)$. If $c>1+L$ and $\M{J}\succ0$, which implies that $\M{G}_k\succ0$, it holds that
\begin{align}\label{SRATE4}
& 2\mathbb{E}\big(\Phi(\bar{\M{X}}_K)-\Phi(\M{X}^*)+\rho\|\M{D}\bar{\M{X}}_K-\M{d}\|\big)\nonumber\\
&\leq \frac{\sigma^2\ln K}{\mu K}+\frac{\mathcal{C}_3}{K}.
\end{align}
\end{Corollary}
{\it Proof}:
See Appendix E.
\hfill{$\blacksquare$}

Assume that $f(\m{x})=\sum_{i=1}^{m}f_i(\m{x})$. Then, to achieve the convergence with constant stepsize, we give the following assumption according to the unified theory of stochastic gradient provided in \cite{SGE1}.
\begin{Assumption}\label{ASS:VR}
There exist $c_1,c_2,c_4\geq0$, $c_3\in(0,1]$, and a random sequence $\{\sigma_k^2\}_{k\geq1}$ such that, for any $k\geq1$, $\mathbb{E}\M{g}^{k}=\nabla F(\M{X}^{k})$,
$\mathbb{E}\|\M{g}^{k}-\nabla F(\M{X}^*)\|^2\leq 2c_1D_F(\M{X}^{k},\M{X}^*)+c_2 \sigma_{k}^2$, and
$\mathbb{E}\sigma_{k+1}^2\leq(1-c_3)\sigma_{k}^2+2c_4D_F(\M{X}^{k},\M{X}^*)$.
\end{Assumption}
This assumption is also used to prove the sublinear convergence rate of the stochastic PDDY (S-PDDY) and stochastic PD3O (S-PD3O) in \cite{DYS-T2}. In terms of Tables 1 and 2 in \cite{SGE1}, Assumption \ref{ASS:VR} is satisfied by several variance-reduced stochastic gradient estimators used in machine learning, such as SAGA \cite{SAGA} and loopless SVRG \cite{LSVRG1,LSVRG2}. Moreover, the value of $c_1,c_2,c_3,c_4$ for corresponding variance-reduced stochastic gradient estimators can be found in \cite[Table 2]{SGE1}. By this assumption, we present the following theorem.
\begin{Theorem}\label{TH3}
Suppose that Assumptions \ref{ASS1} and \ref{ASS:VR} hold. Let $\kappa=\frac{c_2}{c_3}$, $\M{X}_k=\frac{1}{K}\sum_{k=0}^{K-1}\M{X}^{k+1}$, $\tilde{\M{X}}_K=\frac{1}{K}\sum_{k=0}^{K-1}\bar{\M{X}}^{k+1}$, and $\M{X}^*\in\mathcal{X}^*$. If $\alpha\in(0,\frac{1}{2(c_1+\kappa c_4)})$ and $\M{J}\succ 0$, we have
\begin{align}\label{VR:RATE}
&\mathbb{E}(D_F(\M{X}_{K},\M{X}^*)+D_{R}(\tilde{\M{X}}_K,\M{X}^*))\nonumber\\
&\leq \frac{\|\bar{\M{X}}^{0}-(\alpha\M{D}\tr\M{\Lambda}^*+\M{X}^*)\|^2+\|\M{\Lambda}^*\|_{\alpha\M{J}}^2+\kappa\alpha^2\sigma_1^2}{\alpha K},
\end{align}
where $D_F(\M{X}_{K},\M{X}^*)$ is the Bregman divergence of $F(\M{X})$ between $\M{X}_{K}$ and $\M{X}^*$ defined as $D_F(\M{X}_{K},\M{X}^*):=F(\M{X}_{K})-F(\M{X}^*)-\langle \nabla F(\M{X}^*),\M{X}_{K}-\M{X}^*\rangle$, and $D_R(\tilde{\M{X}}_{K},\M{X}^*)$ is the Bregman divergence of $R(\M{X})$ between $\tilde{\M{X}}_{K}$ and $\M{X}^*$.
\end{Theorem}
{\it Proof}:
See Appendix F.
\hfill{$\blacksquare$}

Recall the problem \eqref{ProblemNew} and  its corresponding saddle point problem \eqref{ProblemSaddle}. Let $(\M{X}^*,\M{\Lambda}^*)\in\mathcal{W}^*$, and $\delta^*$ be the conjugate function of the indicator function $\delta_{\{\M{d}\}}(\M{X})$ of the set $\{\M{d}\}$. It is not difficult to verify that
$
\mathcal{L}(\M{X},\M{\Lambda}^*)-\mathcal{L}(\M{X}^*,\M{\Lambda})
=D_F(\M{X},\M{X}^*)+D_R(\M{X},\M{X}^*)+D_{\delta^*}(\M{\Lambda},\M{\Lambda}^*)
=D_F(\M{X},\M{X}^*)+D_R(\M{X},\M{X}^*),
$
where the last equality holds due to $D_{\delta^*}(\M{\Lambda},\M{\Lambda}^*)=0$. Therefore, for $\forall \M{W}\in\mathcal{W}$, by the convexity of $F$ and $R$, one has $\mathcal{L}(\M{X},\M{\Lambda}^*)\leq\mathcal{L}(\M{X}^*,\M{\Lambda}^*)\leq\mathcal{L}(\M{X}^*,\M{\Lambda})$. It implies that the primal--dual gap $\mathcal{L}(\M{X},\M{\Lambda}^*)-\mathcal{L}(\M{X}^*,\M{\Lambda})=D_F(\M{X},\M{X}^*)+D_R(\M{X},\M{X}^*)$ is nonnegative, and if $\M{W}\in\mathcal{W}^*$ the gap $D_F(\M{X},\M{X}^*)+D_R(\M{X},\M{X}^*)$ is zero. Note that $\lim_{k\rightarrow \infty}(\M{X}^{k}-\bar{\M{X}}^k)=0$. Thus, the gap $\mathbb{E}(D_F(\M{X}_{K},\M{X}^*)+D_{R}(\tilde{\M{X}}_K,\M{X}^*))$ is a valid measure. Such measure is also used in \cite{DYS-T2} to present the sublinear convergence rate of S-PDDY and S-PD3O with constant stepsize.

\section{Numerical Experiment}\label{SEC:Simulation}
In this section, we conduct two numerical experiment to validate the obtained theoretical results. All the algorithms are implemented in Matlab R2020b in a computer with 3.30 GHz AMD Ryzen 9 5900HS with Radeon Graphics and 16 GB memory.

\subsection{Generalized Lasso Problem With Equality Constraints}
Consider the following generalized Lasso problem with equality constraints
$$\min_{\m{x}} \frac{1}{2m}\sum_{i=1}^{m}\|\m{A}_i\m{x}-\m{a}_i\|^2+\|\m{B}\m{x}\|_1+\delta_0(\m{Dx}-\m{d}),$$
where $\m{A}_i\in \mathbb{R}^{2n\times n}$, $\m{a}_i\in\mathbb{R}^{2n}$, and each element in $\m{A}_i,\m{a}_i,\m{B},\m{D},\m{d}$ is drawn from the normal distribution $N(0,\tilde{\sigma}^2)$, where $\tilde{\sigma}$ is the standard deviation. We set $m=10$, and $p_1=p_2=20$. In this experiment, we solve this problem by BALPA (S-BALPA) and the existing PD-PSAs including PD3O, PDFP, AFBA, C-V, and TriPD, and compare BALPA (S-BALPA) with the existing PD-PSAs. To implement these PD-PSAs efficiently, we take the specific parameter settings as following. For BALPA (S-BALPA), we set $\alpha={m}/({\sum_{i=1}^{m}\|\m{A}_i\tr\m{A}_i\|})$ and $\gamma=1$. For PDFP, AFBA, C-V, and TriPD, we set $\alpha=1/(\beta\|\M{D}\tr\M{D}\|+\frac{1}{m}{\sum_{i=1}^{m}\|\m{A}_i\tr\m{A}_i\|})$, where $\beta=10^{-3},10^{-6}$ when $\|\M{D}\tr\M{D}\|=O(10^3),O(10^6)$, respectively. For PD3O, we set $\alpha=0.8/(\beta\|\M{D}\tr\M{D}\|+\frac{1}{m}{\sum_{i=1}^{m}\|\m{A}_i\tr\m{A}_i\|})$, where $\beta=10^{-3},10^{-6}$ when $\|\M{D}\tr\M{D}\|=O(10^3),O(10^6)$, respectively. Moreover, we choose SAGA stochastic gradient estimate in S-BALPA.
\begin{table}[!h]
\renewcommand\arraystretch{1.3}
\begin{center}
\caption{Comparison of generalized Lasso problem with equality constraints solved by BALPA (S-BALPA) and the existing PD-PSAs.}
\scalebox{0.95 }{
\begin{tabular}{ccccccc}
\hline
\multirow{2}{*}{\textbf{Algorithm}}&  \multicolumn{2}{c}{$n=2000$}& \multicolumn{2}{c}{$n=4000$}&\multicolumn{2}{c}{$n=6000$}\\
\cline{2-7}
~&Case1&Case2&Case1&Case2&Case1&Case2\\
\hline
PD3O          &    403&1505&387&1443&386&1414  \\
PDFP          &      148&518&144&493&138&478  \\
AFBA          &      141&472&134&407&127&369    \\
C-V, TriPD    &    151 &590& 152& 616& 151& 685\\
BALPA         &    \textbf{15}&\textbf{15}&\textbf{17}&\textbf{17}&\textbf{21}&\textbf{21}  \\
S-BALPA      &  \textbf{5}&\textbf{5}&\textbf{5}&\textbf{5}&\textbf{5}&\textbf{5}   \\
\hline
\end{tabular}}
\label{Table-sim-com1}
\end{center}
\end{table}

In Table \ref{Table-sim-com1}, for various values of $n$ and $\|\M{D}\tr\M{D}\|$, the required numbers of epochs are presented, when the relative error $\|\m{x}^k-\m{x}^*\|/\|\m{x}^0-\m{x}^*\|<10^{-6}$ where $\m{x}^*$ is the optimal solution to the considered problem. Case1 and Case2 denote $\|\M{D}\tr\M{D}\|=O(10^3)$ and $O(10^6)$, respectively. As shown in Table \ref{Table-sim-com1}, BALPA and S-BALPA have a significant acceleration compared with these classic PD-PSAs when $\|\M{D}\tr\M{D}\|$ is large. Moreover, $\|\M{D}\tr\M{D}\|$ does not affect the choice of stepsizes and convergence performance of BALPA (S-BALPA), which is consistent with our theoretical results. To further visualize the numerical results, in Fig. 1, we plot the convergence curves for the considered problem solved by these PD-PSAs, which further demonstrate the efficiency of BALPA and S-BALPA.
\begin{figure}[!h]
  \centering
  \includegraphics[width=1\linewidth]{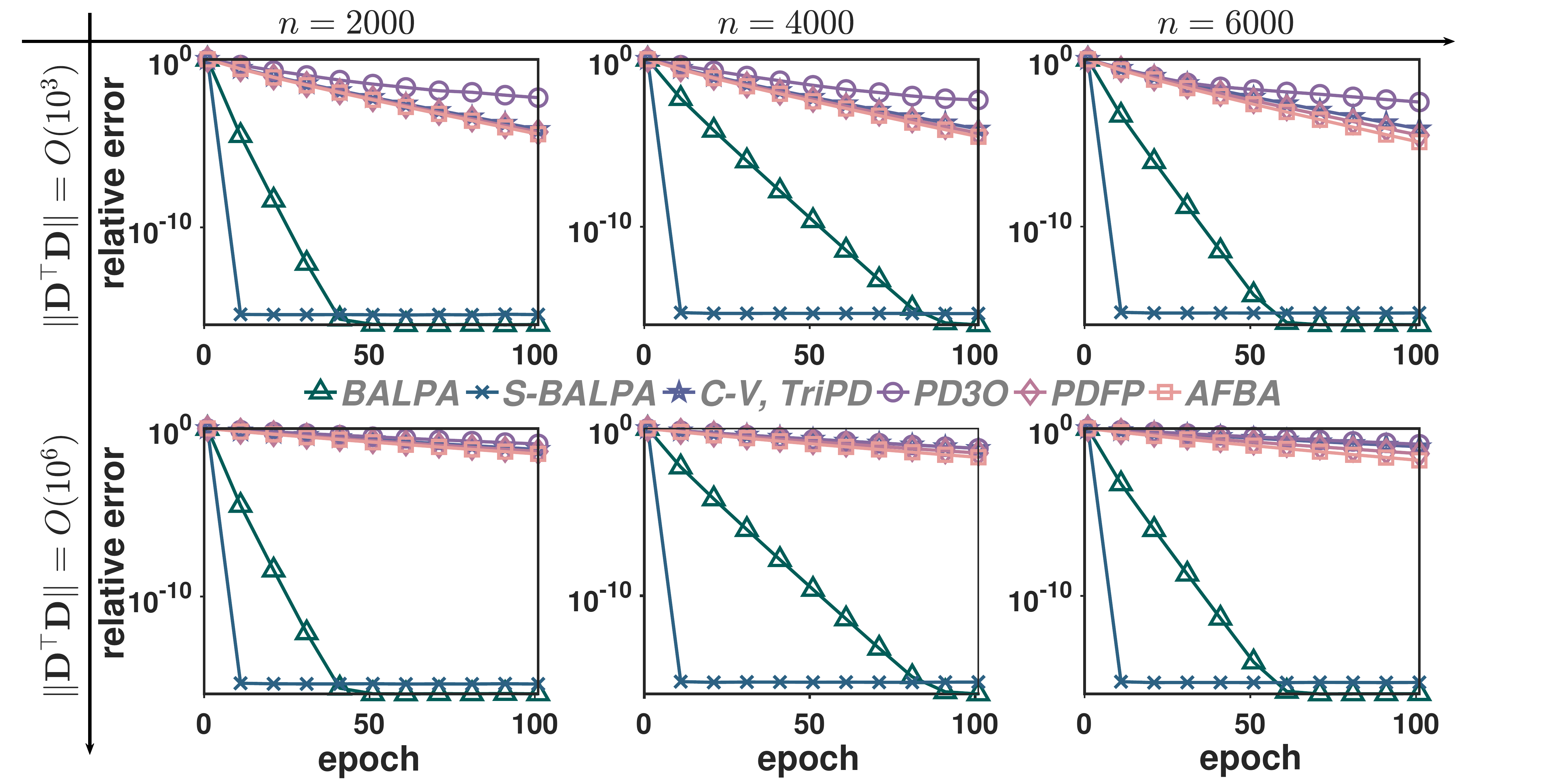}
  \label{Fig-Sim-c1}
  \caption{Convergence curves for generalized Lasso problem with equality constraints solved by BALPA (S-BALPA) and the existing PD-PSAs for various values of $n$ and $\|\M{D}\tr\M{D}\|$.}
\end{figure}
\subsection{Distributed Optimization on Real Datasets}
We use the same setting as \cite{Guo2022}, and consider the distributed logistic regression and linear regression over real datasets. For distributed logistic regression, the optimization problem is as follows
$$
\min _{\tilde{\m{x}}}\sum_{i=1}^{N}\{\frac{1}{m_i}\sum_{j=1}^{m_{i}} \ln (1+e^{-(\mathcal{A}_{i j}^{\top} \tilde{\m{x}})\mathcal{B}_{i j}})+\frac{1}{2}\|\tilde{\m{x}}\|^2+\frac{1}{2}\|\m{B}_i\tilde{\m{x}}\|\}.
$$
For linear regression, the optimization problem is $$
\min _{\tilde{\m{x}}}\sum_{i=1}^{N}\{\frac{1}{m_i}\sum_{j=1}^{m_{i}} \frac{1}{2}\|\mathcal{A}_{ij}\tilde{\m{x}}-\mathcal{B}_{ij}\|^2+\frac{1}{2}\|\tilde{\m{x}}\|^2+\frac{1}{2}\|\m{B}_i\tilde{\m{x}}\|\}.
$$
Here $\tilde{\m{x}}\in \mathbb{R}^{l}$, $\m{B}_i\in \mathbb{R}^{p_1\times l}$ and each element in $\m{B}_i$ is drawn from the normal distribution $N(0,1)$, and any agent $i$ holds its own training date $\left(\mathcal{A}_{i j}, \mathcal{B}_{i j}\right) \in$ $\mathbb{R}^{l} \times\{-1,1\}, j=1, \cdots, m_{i}$, including sample vectors $\mathcal{A}_{i j}$ and corresponding classes $\mathcal{B}_{i j}$. We set the communication topology as a circular graph with 10 nodes, i.e., the 10 agents form a cycle, and use four real datasets including a6a, a9a, covtype, and ijcnn1 \cite{LIBSVM}.

In this experiment, we solve the distributed logistic and linear regression on real datasets by BALPA-Dist, S-BALPA-Dist, DISA, TriPD, PDFP-Dist, and TPUS. We set the primal stepsize $\alpha=0.25$ for these algorithms. To BALPA-Dist, S-BALPA, and DISA, we set $\gamma=\beta=0.5$. To TriPD-Dist, PDFP-Dist, and TPUS, we set $\beta=0.01$. The performance is evaluated by the relative errors $\|\m{x}^k-\m{x}^*\|/\|\m{x}^0-\m{x}^*\|$. The results are illustrated in Fig. 2. It is shown that BALPA-Dist and DISA outperform TriPD-Dist, PDFP-Dist, and TPUS in all the four datasets. The performance for BALPA-Dist and DISA is very close. In this experiment the proximal mapping of the regularizer $\|\tilde{\m{x}}\|$ is easy to compute, and the additional cost of one proximal mapping in DISA can be ignored. Thus, the computational cost for each iteration is also very close for BALPA-Dist and DISA. Moreover, we choose SAGA stochastic gradient estimate in S-BALPA-Dist, and the advantage of using the variance-reduced stochastic gradient estimator is clear.
\begin{figure}[!h]
  \centering
  \includegraphics[width=1\linewidth]{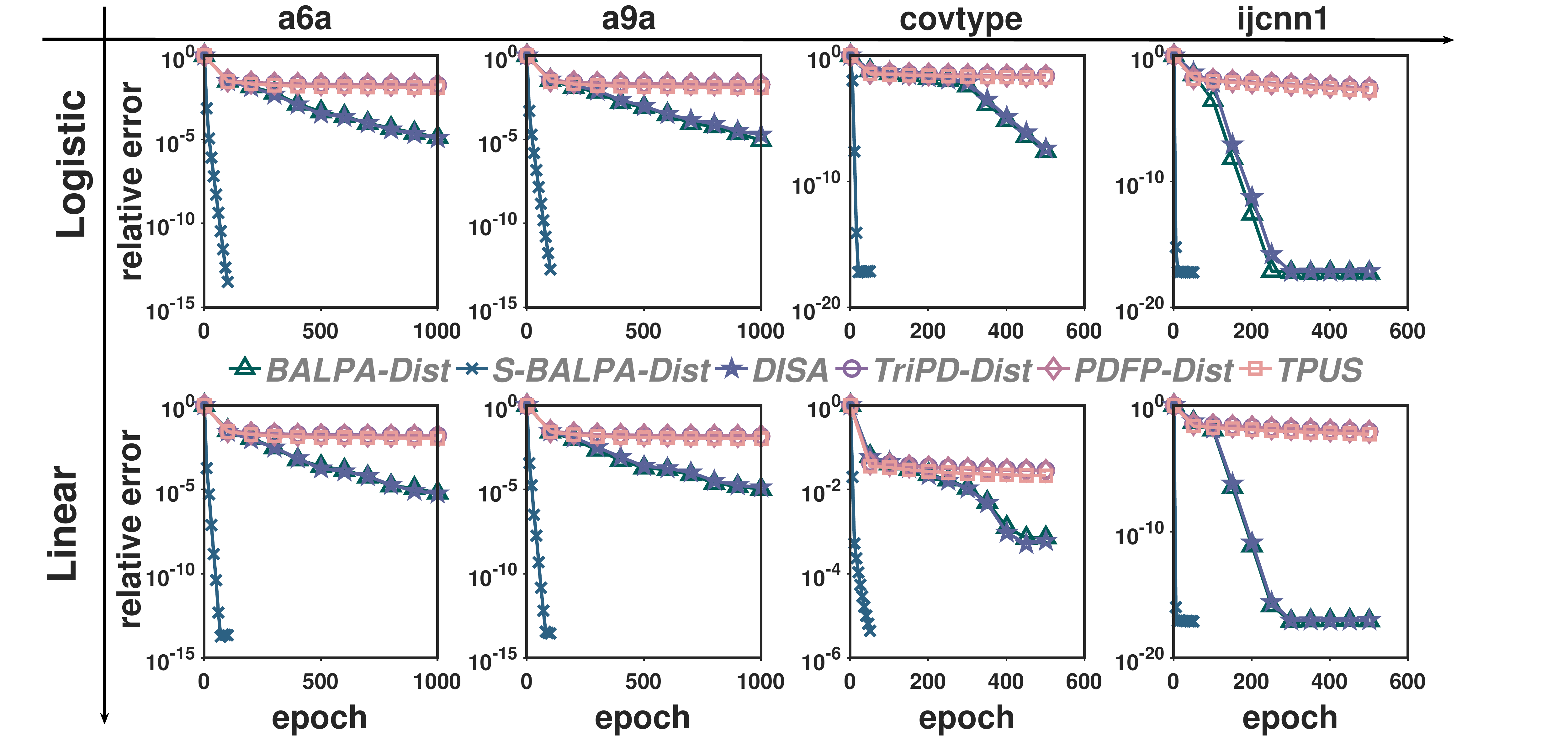}
  \label{Fig-Sim-c2}
  \caption{The results for distributed logistic regression and linear regression over real datasets.}
\end{figure}
\section{Conclusion}\label{SEC:Conclusion}
We have considered a composite optimization problem with linear equality constraints, where the loss function consists of a smooth function and a composition of a nonsmooth proximal function with a linear operator $\m{B}$. For this problem, a novel PD-PSA---named by BALPA has been proposed, which not only overcomes the dependence of the convergence condition on the Euclidean norm of the linear mapping and the equality constraint mapping, but also retains the proximity-induced feature of the classic PD-PSAs. It implies that, unlike existing PD-PSAs, the convergence of BALPA is fast even when these norms are large, and the implementation of BALPA is as simple as existing PD-PSAs. Moreover, we have applied the proposed algorithm to distributed optimization, and developed a new distributed optimization algorithm with better performance than DISA \cite{Guo2022}. In addition, the convergence analysis of the proposed algorithms have been carried out and their superiority have been demonstrated through two numerical experiments. In the future, it would be worthy of investigation to explore the Nesterov- or Newton-type accelerated version of BALPA to improve the convergence rate.

\section*{Appendix A: Proof of Lemma \ref{LEMADA1}}\label{APLEM1}
{\it Proof}:
By \eqref{BAS3} and the definition of proximal mapping, we have that
\begin{align}\label{LEMA1.0}
&R(\M{X})-R(\bar{\M{X}}^k)+\langle \M{X}-\bar{\M{X}}^k,\M{D}\tr\bm{\Lambda}^k+\nabla F(\M{X}^k)+\bm{\delta}^k\nonumber\\
&\quad+\frac{1}{\alpha_k}(\bar{\M{X}}^k-\M{X}^k)\rangle\geq 0, \forall \M{X}\in \mathcal{X}.
\end{align}
Rearranging \eqref{LEMA1.0}, we have that for $\forall \M{X}\in \mathcal{X}$
\begin{align}\label{LEMA1.1}
&\langle\M{X}-\bar{\M{X}}^k,\M{D}\tr(\bm{\Lambda}^k-\bm{\Lambda})+\bm{\delta}^k+\frac{1}{\alpha_k}(\bar{\M{X}}^k-\M{X}^k)\rangle\nonumber\\
& \geq R(\bar{\M{X}}^k)-R(\M{X})+\langle \bar{\M{X}}^k-\M{X},\M{D}\tr\bm{\Lambda}+\nabla F(\M{X}^k)\rangle.
\end{align}
Applying the identity $\langle \m{p}_1-\m{p}_2, \mathrm{H}(\m{q}_1-\m{q}_2)\rangle=\frac{1}{2}\{\|\m{p}_1-\m{q}_2\|^2_{\mathrm{H}}-\|\m{p}_1-\m{q}_1\|^2_{\mathrm{H}}\}+\frac{1}{2}\{\|\m{p}_2-\m{q}_1\|^2_{\mathrm{H}}-\|\m{p}_2-\m{q}_2\|^2_{\mathrm{H}}\}$,
one gets that
\begin{align}\label{LEMA1.2}
&\frac{1}{\alpha_k}\langle\M{X}-\bar{\M{X}}^k,\bar{\M{X}}^k-\M{X}^k\rangle=\frac{1}{2\alpha_k}\|\M{X}-\M{X}^k\|^2\nonumber\\
&\quad-\frac{1}{2\alpha_k}(\|\M{X}-\bar{\M{X}}^k\|^2+\|\bar{\M{X}}^k-\M{X}^k\|^2).
\end{align}
Moreover, it follows from the definition of $\mathbf{1}(\bm{\delta}^k)$ that
\begin{align}\label{LEMA1.3}
&\langle \M{X}-\bar{\M{X}}^k,\bm{\delta}^k \rangle \leq \langle \M{X}-\M{X}^k,\bm{\delta}^k\rangle+\frac{1}{2\eta_k}\|\bm{\delta}^k\|^2\nonumber\\
&\quad+\frac{\mathbf{1}(\bm{\delta}^k)\eta_k}{2}\|\bar{\M{X}}^k-\M{X}^k\|^2.
\end{align}
Thus, combining \eqref{LEMA1.1} with \eqref{LEMA1.2} and \eqref{LEMA1.3}, it gives that
\begin{align}\label{LEMA2}
&\langle\M{X}-\bar{\M{X}}^k,\M{D}\tr(\bm{\Lambda}^k-\bm{\Lambda})\rangle+\langle \M{X}-\M{X}^k,\bm{\delta}^k\rangle\nonumber\\
&\quad+\frac{1}{2\alpha_k}(\|\M{X}-\M{X}^k\|^2-\|\M{X}-\bar{\M{X}}^k\|^2-\|\bar{\M{X}}^k-\M{X}^k\|^2)\nonumber\\
&\quad+\frac{1}{2\eta_k}\|\bm{\delta}^k\|^2+\frac{\mathbf{1}(\bm{\delta}^k)\eta_k}{2}\|\bar{\M{X}}^k-\M{X}^k\|^2\nonumber\\
& \geq R(\bar{\M{X}}^k)-R(\M{X})+\langle \bar{\M{X}}^k-\M{X},\M{D}\tr\bm{\Lambda}+\nabla F(\M{X}^k)\rangle.
\end{align}
It follows from \eqref{BAS1} that
\begin{align}\label{LEMA3}
&\langle\bm{\Lambda}-\bm{\Lambda}^{k+1},-\M{D}(\bar{\M{X}}^{k}-\M{X})+\M{Q}(\bm{\Lambda}^{k+1}-\bm{\Lambda}^k)\rangle\nonumber\\
&\geq\langle\bm{\Lambda}^{k+1}-\bm{\Lambda},\M{d}-\M{D}\M{X}\rangle,\forall \bm{\Lambda}\in \mathcal{Y}.
\end{align}
Similar as \eqref{LEMA1.2}, it holds that
\begin{align*}
&\langle\bm{\Lambda}-\bm{\Lambda}^{k+1},\M{Q}(\bm{\Lambda}^{k+1}-\bm{\Lambda}^k)\rangle=\frac{1}{2}\|\bm{\Lambda}-\bm{\Lambda}^k\|_{\M{Q}}^2\nonumber\\
&\quad-\frac{1}{2}(\|\bm{\Lambda}-\bm{\Lambda}^{k+1}\|_{\M{Q}}^2+\|\bm{\Lambda}^{k+1}-\bm{\Lambda}^k\|_{\M{Q}}^2).
\end{align*}
Summing \eqref{LEMA2} and \eqref{LEMA3}, we have that
\begin{align}\label{LEMA33}
&\langle\M{X}-\bar{\M{X}}^k,\M{D}\tr(\bm{\Lambda}^k-\bm{\Lambda}^{k+1})\rangle+\langle \M{X}-\M{X}^k,\bm{\delta}^k\rangle\nonumber\\
&\quad+\frac{1}{2}(\|\bm{\Lambda}-\bm{\Lambda}^k\|_{\M{Q}}^2-\|\bm{\Lambda}-\bm{\Lambda}^{k+1}\|_{\M{Q}}^2-\|\bm{\Lambda}^{k+1}-\bm{\Lambda}^k\|_{\M{Q}}^2)\nonumber\\
&\quad+\frac{1}{2\alpha_k}(\|\M{X}-\M{X}^k\|^2-\|\M{X}-\bar{\M{X}}^k\|^2-\|\bar{\M{X}}^k-\M{X}^k\|^2)\nonumber\\
&\quad+\frac{1}{2\eta_k}\|\bm{\delta}^k\|^2+\frac{\mathbf{1}(\bm{\delta}^k)\eta_k}{2}\|\bar{\M{X}}^k-\M{X}^k\|^2\nonumber\\
& \geq R(\bar{\M{X}}^k)-R(\M{X})+\langle \bar{\M{X}}^k-\M{X},\M{D}\tr\bm{\Lambda}+\nabla F(\M{X}^k)\rangle\nonumber\\
&\quad+\langle\bm{\Lambda}^{k+1}-\bm{\Lambda},\M{d}-\M{D}\M{X}\rangle, \forall(\M{X},\bm{\Lambda})\in \mathcal{X}\times\mathcal{Y}.
\end{align}
Since $\M{X}^{k+1}=\bar{\M{X}}^{k}+\alpha_k\M{D}\tr(\bm{\Lambda}^k-\bm{\Lambda}^{k+1})$, it gives that
\begin{align*}
&\|\bar{\M{X}}^k-\M{X}\|^2-2\alpha_k\langle\M{X}-\bar{\M{X}}^k,\M{D}\tr(\bm{\Lambda}^k-\bm{\Lambda}^{k+1})\rangle\\
&=\|\bar{\M{X}}^k-\M{X}-\alpha_k\M{D}\tr(\M{\Lambda}^{k+1}-\M{\Lambda}^k)\|^2\\
&\quad-\alpha_k^2\|\M{D}\tr(\M{\Lambda}^k-\M{\Lambda}^{k+1})\|^2\\
&=\|\M{X}^{k+1}-\M{X}\|^2-\alpha_k^2\|\M{D}\tr(\M{\Lambda}^k-\M{\Lambda}^{k+1})\|^2.
\end{align*}
Substituting the above equality into \eqref{LEMA33}, we have
\begin{align}\label{LEMAFF}
&\langle \M{X}-\M{X}^k,\bm{\delta}^k\rangle+\frac{1}{2\eta_k}\|\bm{\delta}^k\|^2-\frac{1}{2}\|\M{\Lambda}^{k+1}-\M{\Lambda}^k\|^2_{\M{Q}-\alpha_k\M{DD}\tr}\nonumber\\
&\quad-\frac{1}{2}(\frac{1}{\alpha_k}-\mathbf{1}(\bm{\delta}^k)\eta_k)\|\bar{\M{X}}^{k}-\M{X}^k\|^2\nonumber\\
&\quad+\frac{1}{2}(\|\bm{\Lambda}-\bm{\Lambda}^k\|_{\M{Q}}^2-\|\bm{\Lambda}-\bm{\Lambda}^{k+1}\|_{\M{Q}}^2)\nonumber\\
&\quad+\frac{1}{2\alpha_k}(\|\M{X}-\M{X}^k\|^2-\|\M{X}-{\M{X}}^{k+1}\|^2)\nonumber\\
& \geq R(\bar{\M{X}}^k)-R(\M{X})+\langle \bar{\M{X}}^k-\M{X},\M{D}\tr\bm{\Lambda}+\nabla F(\M{X}^k)\rangle\nonumber\\
&\quad+\langle\bm{\Lambda}^{k+1}-\bm{\Lambda},\M{d}-\M{D}\M{X}\rangle, \forall(\M{X},\bm{\Lambda})\in \mathcal{X}\times\mathcal{Y}.
\end{align}
Since $F$ is convex and $L$-smooth, by \eqref{SM1}, it holds that
\begin{align}\label{L2SMOOTH}
&\big\langle \M{\bar{X}}^{k}-\M{X},\nabla F(\M{X}^k) \big\rangle\nonumber\\
&=\big\langle \M{\bar{X}}^{k}-\M{X},\nabla F(\M{X}^k)- \nabla F(\M{X})\big\rangle+\big\langle \M{\bar{X}}^{k}-\M{X},\nabla F(\M{X})\big\rangle\nonumber\\
&\geq-\frac{L}{4} \|\bar{\M{X}}^{k}-\M{X}^k\|^2+\big\langle \M{\bar{X}}^{k}-\M{X},\nabla F(\M{X})\big\rangle.
\end{align}
Using the notations of $\M{H}_k$ and $\M{M}_k$, and combining and \eqref{LEMAFF} with \eqref{L2SMOOTH}, we get the assertion \eqref{FunA1} immediately.

Moreover, since $f$ is convex and $L$-smooth, by \eqref{SM2}, one has
\begin{align*}
\big\langle \M{X}-\bar{\M{X}}^{k},\nabla F(\M{X}^k)\rangle\leq F(\M{X})-F(\bar{\M{X}}^{k})+\frac{L}{2}\|\M{X}^k-\bar{\M{X}}^{k}\|^2.
\end{align*}
It follows from \eqref{LEMAFF} and notations of $\M{H}_k$ and $\M{G}_k$ that the assertion \eqref{FunA2} holds.
\hfill{$\blacksquare$}

\section*{Appendix B: Proof of Theorem \ref{TH1}}\label{APTH1}
{\it Proof}:
Since $\sigma=0$, we have $\tau_k\equiv\frac{1}{\alpha}$, $\M{H}_k\equiv\M{H}$, $\M{M}_k\equiv \M{M}$, and $\M{G}_k\equiv \M{G}$, where $\M{H}: \mathcal{W}\rightarrow \mathcal{W}:(\M{X},\M{\Lambda})\mapsto (\frac{1}{\alpha}\M{X},\M{Q}\M{\Lambda})$,
\begin{align*}
\M{M}: \mathcal{W}\rightarrow \mathcal{W}:&(\M{X},\M{\Lambda})\mapsto ((\frac{1}{\alpha}-\frac{L}{2})\M{X},(\M{Q}-\alpha\M{DD}\tr)\M{\Lambda}),\\
\M{G}: \mathcal{W}\rightarrow \mathcal{W}:&(\M{X},\M{\Lambda})\mapsto ((\frac{1}{\alpha}-L)\M{X},(\M{Q}-\alpha\M{DD}\tr)\M{\Lambda}),
\end{align*}
are self-adjoint linear mappings. For any $\M{W}^*\in\mathcal{W}^*$, letting $\M{W}=\M{W}^*$ in \eqref{FunA1}, it holds that
\begin{align*}
&0\overset{\eqref{IV2}}{\leq}2R(\bar{\M{X}}^k)-2R(\M{X}^*)+2\langle\mathbf{V}^{k+1}-\mathbf{W}^*,(\bm{B}+\bm{C})(\mathbf{W}^*)\rangle\\
&\leq \|\M{W}^*-\M{W}^k\|^2_{\M{H}}-\|\M{W}^*-\M{W}^{k+1}\|^2_{\M{H}}-\|\M{V}^{k+1}-\M{W}^k\|^2_{\M{M}},
\end{align*}
which implies that the sequence $\{\M{W}^k\}_{k\geq0}$ generated by BALPA is a Fej\'{e}r monotone sequence with respect to $\mathcal{W}^*$ in $\mathbf{H}$-norm. Since $\M{H}\succ0$, for any $\M{W}^*\in\mathcal{W}^*$, the sequence $\{\|\M{W}^k-\M{W}^*\|\}_{k\geq0}$ is bounded.

Summing the above inequality over $k=0,1,\cdots,K-1$, it gives that
\begin{align*}
&\sum_{k=0}^{K-1}\|\M{W}^k-\M{V}^{k+1}\|^2_{\M{M}}\\
&\leq \sum_{k=0}^{K-1}(\|\M{W}^k-\M{W}\|_{\M{H}}^2-\|\M{W}^{k+1}-\M{W}\|_{\M{H}}^2)\\
&=\|\M{W}^0-\M{W}^*\|_{\M{H}}^2-\|\M{W}^{K}-\M{W}^*\|_{\M{H}}^2,\quad \forall K\geq1.
\end{align*}
Thus, it holds that for any $K\geq1$,
\begin{align}\label{THPPPFF1}
\sum_{k=0}^{K-1}\|\M{W}^k-\M{V}^{k+1}\|^2_{\M{M}}\leq \|\M{W}^0-\M{W}^*\|_{\M{H}}^2,
\end{align}
which implies that $\sum_{k=0}^{\infty}\|\M{W}^k-\M{V}^{k+1}\|^2_{\M{M}}<\infty$. Since $\M{M}$ is positive definite, we have $\lim_{k\rightarrow \infty}\|\M{W}^k-\M{V}^{k+1}\|=0$. Let $\mathbf{W}^{\infty}=[(\M{X}^{\infty})\tr,(\M{\Lambda}^{\infty})]\tr$ be an accumulation point of $\{\mathbf{W}^k\}$ and $\{\mathbf{W}^{k_j}\}$ be a subsequence converging
to $\mathbf{W}^{\infty}$. By \eqref{BAS2}, we have
$\|\bar{\M{X}}^{k}-\M{X}^{k+1}\|=\|\alpha\M{D}\tr(\M{\Lambda}^k-\M{\Lambda}^{k+1})\|$. Hence, it holds that
$
\lim_{k\rightarrow \infty}\|\M{V}^{k+1}-\M{W}^{k+1}\|=\lim_{k\rightarrow \infty} (\|\bar{\M{X}}^{k}-\M{X}^{k+1}\|+\|\M{\Lambda}^{k+1}-\M{\Lambda}^{k+1}\|)=0,
$
which implies that $\mathbf{w}^{\infty}$ is also an accumulation point of $\{\mathbf{V}^k\}$ and $\{\mathbf{V}^{k_j}\}$ is a subsequence converging
to $\mathbf{W}^{\infty}$. By \eqref{LEMA1.0} and \eqref{LEMA3}, one has that for any $\M{W}\in\mathcal{W}$
\begin{align*}
&R(\M{X})-R(\bar{\M{X}}^{k_j})+\langle \M{X}-\bar{\M{X}}^{k_j},\M{D}\tr\bm{\Lambda}^{k_j}+\nabla F(\M{X}^{k_j})\\
&\quad+\frac{1}{\alpha}(\bar{\M{X}}^{k_j}-\M{X}^{k_j})\rangle\geq 0, \\
&\big\langle \M{\Lambda}-\M{\Lambda}^{k_j+1},\M{d}-\M{D}\bar{\M{X}}^{k_j}+\M{Q}(\M{\Lambda}^{k_j+1}-\M{\Lambda}^{k_j})\big\rangle\geq0.
\end{align*}
Taking $k_j\rightarrow\infty$ in the above two inequalities, it holds that for any $\M{W}\in\mathcal{W}$
\begin{align*}
R(\M{X})-R(\M{X}^{\infty})+\big\langle\mathbf{W}-\mathbf{W}^{\infty},(\bm{B}+\bm{C})(\mathbf{W}^{\infty})\big\rangle\geq 0.
\end{align*}
Compared with \eqref{IV2}, we can deduce that $\M{W}^{\infty} \in \mathcal{W}^*$. Therefore, by the Fej\'{e}r monotonicity of $\{\M{W}^k\}_{k\geq0}$, we have that $\|\M{W}^k-\M{W}^{\infty}\|_{\mathbf{H}}^2-\|\M{W}^{k+1}-\M{W}^{\infty}\|_{\mathbf{H}}^2\geq 0$.
It implies that the sequence $\{\|\M{W}^k-\M{W}^{\infty}\|_{\mathbf{H}}\}_{k\geq0}$ converges to a unique limit point. Then, with $\M{W}^{\infty}$ being an accumulation point of $\{\M{W}^k\}$, we have $\lim_{k\rightarrow \infty}\M{W}^k=\M{W}^{\infty}$.
\hfill{$\blacksquare$}

\section*{Appendix C: Proof of Theorem \ref{TH:RATE1}}
{\it Proof}: It follows from \eqref{THPPPFF1} that
\begin{align*}
&(K+1)\min_{0\leq k\leq K}\{\|\M{V}^{k+1}-\M{W}^k\|^2_{\M{M}}\}\leq\sum_{k=0}^{K}\|\M{W}^k-\M{V}^{k+1}\|^2_{\M{M}}\\
&\leq\sum_{k=0}^{\infty}\|\M{W}^k-\M{V}^{k+1}\|^2_{\M{M}}\leq \|\M{W}^0-\M{W}^*\|_{\M{H}}^2<\infty,
\end{align*}
which implies that the assertion \eqref{D-BPDA-RATA1} holds. It follows from \eqref{FunA2} that for $\forall \M{W}\in \mathcal{W}$
\begin{align*}
&2\mathcal{L}(\bar{\M{X}}^k,\M{\Lambda})-2\mathcal{L}(\M{X},\M{\Lambda}^{k+1})\nonumber\\
&\leq\|\M{W}-\M{W}^k\|^2_{\M{H}}-\|\M{W}-\M{W}^{k+1}\|^2_{\M{H}}-\|\M{V}^{k+1}-\M{W}^k\|^2_{\M{G}}.
\end{align*}
Since $\alpha\in(0,\frac{1}{L})$, the linear mapping $\M{G}\succ0$. Then, summing the above inequality over $k=0,1,\cdots,K-1$, we obtain
$
2\sum_{k=0}^{K-1}(\mathcal{L}(\bar{\M{X}}^k,\M{\Lambda})-\mathcal{L}(\M{X},\M{\Lambda}^{k+1}))\leq\|\M{W}^0-\M{W}\|_{\M{H}}^2-\|\M{W}^K-\M{W}\|_{\M{H}}^2.
$
By the convexity of $F$, $R$ and the definition of $\bar{\M{X}}_K=\frac{1}{K}\sum_{k=0}^{K-1}\bar{\M{X}}^k$, ${\M{\Lambda}}_K=\frac{1}{K}\sum_{k=0}^{K-1}{\M{\Lambda}}^{k+1}$, we have
$$
K(\mathcal{L}(\bar{\M{X}}_K,\M{\Lambda})-\mathcal{L}(\M{X},{\M{\Lambda}}_K))\leq\sum_{k=0}^{K-1}(\mathcal{L}(\bar{\M{X}}^k,\M{\Lambda})-\mathcal{L}(\M{X},\M{\Lambda}^{k+1})).
$$
Therefore, the primal--dual gap \eqref{PDGAP} holds.

Note the inequality \eqref{PDGAP} is true for all $\M{X}\in\mathcal{X}$, $\M{\Lambda}\in \mathcal{Y}$, hence it is also true for any optimal solution $\M{X}^*$ and $\mathcal{B}_{\rho}=\{\M{\Lambda}:\|\M{\Lambda}\|\leq\rho\}$, where $\rho>0$ is any given positive number. Note that the mapping $\bm{B}$ is affine with a skew-symmetric matrix and thus we have $\langle \M{W}_1-\M{W}_2,\bm{B}\M{W}_1-\bm{B}\M{W}_2 \rangle\equiv0,\forall \M{W}_1,\M{W}_2\in \mathcal{W}$. As a result, by letting $\M{W}=[(\M{X}^*)\tr,\M{\Lambda}\tr]\tr$ and $\M{V}_K=[(\bar{\M{X}}_K)\tr,(\M{\Lambda}_K)\tr]\tr$, it follows that
\begin{align}\label{THMproof17}
&\sup_{\M{\Lambda}\in\mathcal{B}_{\rho}}\big\{\mathcal{L}(\bar{\M{X}}_K,\M{\Lambda})-\mathcal{L}(\M{X}^*,{\M{\Lambda}}_K)\big\}\nonumber\\
&=\sup_{\M{\Lambda}\in\mathcal{B}_{\rho}}\big\{\Phi(\bar{\M{X}}_K)-\Phi(\M{X}^*)+\langle{\mathbf{V}}_{K}-\mathbf{W},\bm{B}\M{W}\rangle\big\}\nonumber\\
&=\sup_{\M{\Lambda}\in\mathcal{B}_{\rho}}\big\{\Phi(\bar{\M{X}}_K)-\Phi(\M{X}^*)+\langle{\mathbf{V}}_{K}-\mathbf{W},\bm{B}{\mathbf{V}}_{K}\rangle\big\}\nonumber\\
&=\sup_{\M{\Lambda}\in\mathcal{B}_{\rho}}\big\{\Phi(\bar{\M{X}}_K)-\Phi(\M{X}^*)+\langle\bar{\M{X}}_K-\M{X}^*,\M{D}\tr{\M{\Lambda}}_K\rangle\nonumber\\
&\quad\quad\quad\quad+\langle {\M{\Lambda}}_K-\M{\Lambda},\M{d}-\M{D}\bar{\M{X}}_K\rangle\big\}\nonumber\\
&=\sup_{\M{\Lambda}\in\mathcal{B}_{\rho}}\big\{\Phi(\bar{\M{X}}_K)-\Phi(\M{X}^*)+\langle \M{\Lambda},\M{D}\tilde{\M{X}}_K-\M{d}\rangle\big\}\nonumber\\
&=\Phi(\bar{\M{X}}_K)-\Phi(\M{X}^*)+\rho\|\M{D}\bar{\M{X}}_K-\M{d}\|.
\end{align}
Note that $\M{\Lambda}^0=0$. Combining \eqref{PDGAP} and \eqref{THMproof17}, we have
\begin{align}\label{THMproof18}
&\Phi(\bar{\M{X}}_K)-\Phi(\M{X}^*)+\rho\|\M{D}\bar{\M{X}}_K-\M{d}\|\nonumber\\
&\leq\frac{1}{2K}\sum_{k=0}^{K-1}\frac{1}{\alpha}(\|\M{X}^k-\M{X}^*\|^2-\|\M{X}^{k+1}-\M{X}^*\|^2)\nonumber\\
&\quad+\frac{1}{2K}\sup_{\M{\Lambda}\in\mathcal{B}_{\rho}}\big\{\sum_{k=0}^{K-1}(\|\M{\Lambda}^k-\M{\Lambda}\|^2_{\M{Q}}-\|\M{\Lambda}^{k+1}-\M{\Lambda}\|^2_{\M{Q}})\big\}\nonumber\\
&=\frac{1}{2K}(\frac{1}{\alpha}\|\M{X}^0-\M{X}^*\|-\|\M{X}^{K}-\M{X}^*\|)\nonumber\\
&\quad+\frac{1}{2K}\sup_{\M{\Lambda}\in\mathcal{B}_{\rho}}\big\{\|\M{\Lambda}^0-\M{\Lambda}^*\|^2_{\M{Q}}-\|\M{\Lambda}^K-\M{\Lambda}\|^2_{\M{Q}}\big\}\nonumber\\
&\leq \frac{\frac{1}{\alpha}\|\M{X}^0-\M{X}^*\|+\rho^2\|\M{Q}\|}{2K}.
\end{align}
Let $\omega_{K}\triangleq\frac{\frac{1}{\alpha}\|\M{X}^0-\M{X}^*\|+\rho^2\|\M{Q}\|}{2K}$.
For any $\rho>\|\M{\Lambda}^*\|$, applying \cite[Lemma 2.3]{ADMM-D1} in \eqref{THMproof18}, it gives that
\begin{align*}
\|\M{D}\bar{\M{X}}_K-\M{d}\|&\leq\frac{\omega_{K}}{\rho-\|\M{\Lambda}^*\|},\\
-\frac{\omega_{K}\|\M{\Lambda}^*\|}{\rho-\|\M{\Lambda}^*\|}&\leq\Phi(\bar{\M{X}}_K)-\Phi(\M{X}^*)\leq \omega_{K}.
\end{align*}
Letting $\rho=\max\{1+\|\M{\Lambda}^*\|,2\|\M{\Lambda}^*\|\}$, \eqref{D-BPDA-RATA2} holds.
\hfill{$\blacksquare$}

\section*{Appendix D: Proof of Theorem \ref{TH2}}\label{APTH2}
{\it Proof}:
Letting $\M{W}=\M{W}^*$ in \eqref{FunA1}, it holds that
\begin{align*}
0\overset{\eqref{IV2}}{\leq}&2R(\bar{\M{X}}^k)-2R(\M{X}^*)+2\langle\mathbf{V}^{k+1}-\mathbf{W}^*,(\bm{B}+\bm{C})\mathbf{W}^*\rangle\\
\leq& 2\langle \M{X}-\M{X}^k,\bm{\delta}^k\rangle+\frac{1}{\eta_k}\|\bm{\delta}^k\|^2-\|\M{V}^{k+1}-\M{W}^k\|^2_{\M{M}_k}\\
&+\|\M{W}^*-\M{W}^k\|^2_{\M{H}_k}-\|\M{W}^*-\M{W}^{k+1}\|^2_{\M{H}_k}.
\end{align*}
Summing the above inequality over $k=0,1,\cdots,K-1$, it deduces that
\begin{align}\label{SRATEPROOF11}
&\sum_{k=0}^{K-1}\|\M{V}^{k+1}-\M{W}^k\|^2_{\M{M}_k}\leq\sum_{k=0}^{K-1}(2\langle \M{X}-\M{X}^k,\bm{\delta}^k\rangle+\frac{\|\bm{\delta}^k\|^2}{\eta_k})\nonumber\\
&+\sum_{k=0}^{K-1}(\|\M{W}^*-\M{W}^k\|^2_{\M{H}_k}-\|\M{W}^*-\M{W}^{k+1}\|^2_{\M{H}_k}).
\end{align}
When $\mathcal{D}<\infty$, let $\eta_k=1+\sqrt{k}$, $\frac{1}{\alpha_k}=c+\sqrt{k}$. We have
\begin{align*}
&\sum_{k=0}^{K-1}\frac{1}{\alpha_k}(\|\M{X}^k-\M{X}^*\|^2-\|\M{X}^{k+1}-\M{X}^*\|^2)\nonumber\\
&=\frac{1}{\alpha_0}\|\M{X}^0-\M{X}^*\|^2-\frac{1}{\alpha_{K-1}}\|\M{X}^K-\M{X}^*\|^2\nonumber\\
&\quad+\sum_{k=0}^{K-2}(\sqrt{k+1}-\sqrt{k})\|\M{X}^{k+1}-\M{X}^*\|^2\nonumber\\
&\leq c\|\M{X}^0-\M{X}^{*}\|^2+\sqrt{K-1}\mathcal{D}^2,
\end{align*}
When $\mathcal{D}=\infty$, let $\eta_k=1+\sqrt{K-1}$, $\frac{1}{\alpha_k}=c+\sqrt{K-1}$. It holds that
\begin{align*}
&\sum_{k=0}^{K-1}\frac{1}{\alpha_k}(\|\M{X}^k-\M{X}^*\|^2-\|\M{X}^{k+1}-\M{X}^*\|^2)\nonumber\\
&=\sum_{k=0}^{K-1}(c+\sqrt{K-1})(\|\M{X}^k-\M{X}^*\|^2-\|\M{X}^{k+1}-\M{X}^*\|^2)\nonumber\\
&\leq c\|\M{X}^0-\M{X}^{*}\|^2+\sqrt{K-1}\|\M{X}^0-\M{X}^{*}\|^2.
\end{align*}
From the definition of $\bar{\mathcal{D}}$, it holds that
\begin{align}\label{LpPPROOF1}
&\sum_{k=0}^{K-1}\frac{1}{\alpha_k}(\|\M{X}^k-\M{X}^*\|^2-\|\M{X}^{k+1}-\M{X}^*\|^2)\nonumber\\
&\leq c\|\M{X}^0-\M{X}^{*}\|^2+\sqrt{K-1}\bar{\mathcal{D}}^2,
\end{align}
which implies that
\begin{align}\label{SRATEPROOF22}
&\sum_{k=0}^{K-1}\big(\|\M{W}^*-\M{W}^k\|^2_{\M{H}_k}-\|\M{W}^*-\M{W}^{k+1}\|^2_{\M{H}_k}\big)\nonumber\\
&\leq \frac{1}{\alpha_0}\|\M{X}^0-\M{X}^{*}\|^2+\|\M{\Lambda}^*\|^2_{\M{Q}}+\sqrt{K-1}\bar{\mathcal{D}}^2.
\end{align}
Since $\bm{\delta}^k=\M{g}^k-\nabla F(\M{X}^k)$, it gives that
$\mathbb{E}(\bm{\delta}^k)=0, ~\mathbb{E}(\|\bm{\delta}^k\|^2)\leq \sigma^2,\forall k\geq0$.
When $\eta_k=1+\sqrt{k}$, one has
$
\sum_{k=0}^{K-1}\frac{1}{\eta_k}=\sum_{k=0}^{K-1}\frac{1}{1+\sqrt{k}}\leq\int_{0}^{K-1}\frac{1}{1+\sqrt{s}}ds\leq\sqrt{K},
$
and when $\eta_k=1+\sqrt{K-1}$, one has
$\sum_{k=0}^{K-1}\frac{1}{\eta_k}=\sum_{k=0}^{K-1}\frac{1}{1+\sqrt{K-1}}\leq \sqrt{K}.$
One can deduce that
\begin{align}\label{SRATEPROOF33}
\mathbb{E}\Big(\sum_{k=0}^{K-1}\big(2\langle \M{X}-\M{X}^k,\bm{\delta}^k\rangle+\frac{\|\bm{\delta}^k\|^2}{\eta_k}\big)\Big)\leq \sqrt{K}\sigma^2.
\end{align}
Therefore, taking the expectation of both sides in \eqref{SRATEPROOF11}, then using \eqref{SRATEPROOF22} and \eqref{SRATEPROOF33}, one has that
\begin{align*}
&\max_{0\leq k\leq K-1}\big\{\mathbb{E}\big(\|\M{V}^{k+1}-\M{W}^k\|^2_{\M{M}_k}\big)\big\}\\
&\leq\frac{1}{K}\mathbb{E}\Big(\sum_{k=0}^{K-1}\|\M{V}^{k+1}-\M{W}^k\|^2_{\M{M}_k}\Big)\\
&\leq\frac{\bar{\mathcal{D}}^2+\sigma^2}{\sqrt{K}}+\frac{c\|\M{X}^0-\M{X}^{*}\|^2+\|\M{\Lambda}^*\|^2_{\M{Q}}}{K},
\end{align*}
i.e., the assertion \eqref{SRATE1} holds.

%According to the updates of $\bar{\M{X}}^{k}$ and $\M{\Lambda}^{k+1}$, we have
%\begin{align*}
%&0\in \partial R(\bar{\M{X}}^{k})+\nabla F(\M{X}^k)+\bm{\delta}^k+{\M{D}}\tr\M{\Lambda}^k+\frac{1}{\alpha_k}(\bar{\M{X}}^{k}-\M{X}^k),\\
%&0=\M{d}-{\M{D}}\bar{\M{X}}^{k}+\M{Q}(\M{\Lambda}^{k+1}-\M{\Lambda}^{k}),
%\end{align*}
%which implies that
%\begin{align*}
%&\nabla F(\bar{\M{X}}^{k})-\nabla F(\M{X}^{k})+{\M{D}}\tr(\M{\Lambda}^{k+1}-\M{\Lambda}^{k})-\frac{1}{\alpha_k}(\bar{\M{X}}^{k}-\M{X}^k)\\
%&+\bm{\delta}^k\in\partial R(\bar{\M{X}}^{k})+\nabla F(\bar{\M{X}}^{k})+{\M{D}}\tr\M{\Lambda}^{k+1},\\
%&-\M{Q}(\M{\Lambda}^{k+1}-\M{\Lambda}^{k})=\M{d}-{\M{D}}\bar{\M{X}}^{k}.
%\end{align*}
%Therefore
%\begin{align*}
%&\mathrm{dist}^2(0,\mathcal{K}(\M{V}^{k+1}))\leq\|\nabla F(\bar{\M{X}}^{k})-\nabla F(\M{X}^{k})+{\M{D}}\tr(\M{\Lambda}^{k+1}-\M{\Lambda}^{k})\\
%&\quad+\frac{1}{\alpha_k}(\bar{\M{X}}^{k}-\M{X}^k)+\M{\delta}^k\|^2+\|\M{Q}(\M{\Lambda}^{k+1}-\M{\Lambda}^{k})\|^2\\
%&\leq\kappa^2_1\|\M{V}^{k+1}-\M{W}^k\|^2_{\M{M}_k},
%\end{align*}
%where $\kappa^2_1=$

Since $\M{M}_k\succ 0$, summing the inequality \eqref{FunA2} over $k=0,1,\cdots,K-1$, we have that for any $\M{W}\in\mathcal{W}$
\begin{align}\label{LEM22A1}
&\sum_{k=0}^{K-1}( \Phi(\bar{\M{X}}^k)-\Phi(\M{X})+\langle\mathbf{V}^{k+1}-\mathbf{W},\bm{B}\mathbf{W}\rangle)\nonumber\\
&\leq \sum_{k=0}^{K-1}\langle \M{X}-\M{X}^k,\bm{\delta}^k\rangle+\sum_{k=0}^{K-1}\frac{1}{2\eta_k}\|\bm{\delta}^k\|^2\nonumber\\
&+\sum_{k=0}^{K-1}\big(\frac{1}{2}\|\M{W}-\M{W}^k\|^2_{\M{H}_k}-\frac{1}{2}\|\M{W}-\M{W}^{k+1}\|^2_{\M{H}_k}\big).
\end{align}
By the convexity of $\Phi$ and the definition of $\bar{\M{X}}_K$, it holds that $K\Phi(\bar{\M{X}}_K)\leq \sum_{k=0}^{K-1}\Phi(\bar{\M{X}}^k)$. Since the mapping $\bm{B}$ is affine with a skew-symmetric matrix and thus we have $\langle \M{W}_1-\M{W}_2,\bm{B}\M{W}_1-\bm{B}\M{W}_2 \rangle\equiv0,\forall \M{W}_1,\M{W}_2\in \mathcal{W}$. Therefore, it follows from \eqref{LEM22A1} that for any $\M{W}\in \mathcal{W}$
\begin{align}\label{LEM22A2}
&2K( \Phi(\bar{\M{X}}_K)-\Phi(\M{X})+\langle{\mathbf{V}}_{K}-\mathbf{W},\bm{B}{\mathbf{V}}_{K}\rangle)\nonumber\\
&\leq 2\sum_{k=0}^{K-1}\langle \M{X}-\M{X}^k,\bm{\delta}^k\rangle+\sum_{k=0}^{K-1}\frac{1}{\eta_k}\|\bm{\delta}^k\|^2\nonumber\\
&+\sum_{k=0}^{K-1}(\|\M{W}-\M{W}^k\|^2_{\M{H}_k}-\|\M{W}-\M{W}^{k+1}\|^2_{\M{H}_k}).
\end{align}
Combining \eqref{THMproof17} and \eqref{LEM22A2}, we have that for any $\rho>0$
\begin{align}\label{LAMPROOFAB2}
&2K(\Phi(\bar{\M{X}}_K)-\Phi(\M{X}^*)+\rho\|\M{D}\bar{\M{X}}_K-\M{d}\|)\nonumber\\
&\leq 2\sum_{k=0}^{K-1}\langle \M{X}^*-\M{X}^k,\bm{\delta}^k\rangle+\sum_{k=0}^{K-1}\frac{1}{\eta_k}\|\bm{\delta}^k\|^2\nonumber\\
&+\sum_{k=0}^{K-1}\frac{1}{\alpha_k}(\|\M{X}^k-\M{X}^*\|^2-\|\M{X}^{k+1}-\M{X}^*\|^2)\nonumber\\
&+\sup_{\M{\Lambda}\in\mathcal{B}_{\rho}}\big\{\sum_{k=0}^{K-1}(\|\M{\Lambda}^k-\M{\Lambda}\|^2_{\M{Q}}-\|\M{\Lambda}^{k+1}-\M{\Lambda}\|^2_{\M{Q}})\big\}.
\end{align}
Note that
\begin{align}\label{LpPPROOF2}
&\sup_{\M{\Lambda}\in\mathcal{B}_{\rho}}\big\{\sum_{k=0}^{K-1}(\|\M{\Lambda}^k-\M{\Lambda}\|^2_{\M{Q}}-\|\M{\Lambda}^{k+1}-\M{\Lambda}\|^2_{\M{Q}})\big\}\nonumber\\
&=\sup_{\M{\Lambda}\in\mathcal{B}_{\rho}}\big\{\|\M{\Lambda}^0-\M{\Lambda}\|^2_{\M{Q}}-\|\M{\Lambda}^K-\M{\Lambda}\|^2_{\M{Q}}\big\}\nonumber\\
&\leq\sup_{\M{\Lambda}\in\mathcal{B}_{\rho}}\{\|\M{\Lambda}^0-\M{\Lambda}\|^2_{\M{Q}}\}\leq\sup_{\M{\Lambda}\in\mathcal{B}_{\rho}}\{\|\M{Q}\|\|\M{\Lambda}\|^2\big\}=\rho^2\|\M{Q}\|,
\end{align}
where the last inequality follows from $\M{\Lambda}^0=0$. Therefore, by plugging \eqref{LpPPROOF1} and \eqref{LpPPROOF2} into \eqref{LAMPROOFAB2}, it yields
\begin{align}\label{LEM22A3}
& \Phi(\bar{\M{X}}_K)-\Phi(\M{X}^*)+\rho\|\M{D}\bar{\M{X}}_K-\M{d}\|\nonumber\\
&\leq \frac{1}{K}\sum_{k=0}^{K-1}\langle \M{X}^*-\M{X}^k,\bm{\delta}^k\rangle+\frac{1}{K}\sum_{k=0}^{K-1}\frac{1}{2\eta_k}\|\bm{\delta}^k\|^2\nonumber\\
&\quad+\frac{c\|\M{X}^0-\M{X}^{*}\|^2+\sqrt{K-1}\bar{\mathcal{D}}^2+\rho^2\|\M{Q}\|}{2K}.
\end{align}
Since $\bm{\delta}^k=\M{g}^k-\nabla F(\M{X}^k)$, it gives that
$
\mathbb{E}(\bm{\delta}^k)=0, ~\mathbb{E}(\|\bm{\delta}^k\|^2)\leq \sigma^2,\forall k\geq0.
$
Taking the expectation of both sides in \eqref{LEM22A3} and using the inequality $\sum_{k=0}^{K-1}\frac{1}{\eta_k}\leq \sqrt{K}$ ,
we get \eqref{SRATE2} immediately.
\hfill{$\blacksquare$}

\section*{Appendix E: Proof of Corollary \ref{COR:RATE1}}
{\it Proof}:
Since $f$ is $\mu$-strongly convex and $L$-smooth, one has
$\langle \M{X}-\bar{\M{X}}^{k},\nabla F(\M{X}^k)\rangle \leq F(\M{X})-F(\bar{\M{X}}^{k})+\frac{L}{2}\|\M{X}^k-\bar{\M{X}}^{k}\|^2-\frac{\mu}{2}\|\M{X}-\M{X}^k\|^2$.
By \eqref{LEMAFF} and notations of $\M{H}_k$ and $\M{G}_k$, it holds that for $\forall \M{W}\in \mathcal{W}$
\begin{align*}
&\Phi(\bar{\M{X}}^k)-\Phi(\M{X})+\langle\mathbf{V}^{k+1}-\mathbf{W},\bm{B}\mathbf{W}\rangle+\frac{\mu}{2}\|\M{X}-\M{X}^k\|^2\nonumber\\
&\leq \langle \M{X}-\M{X}^k,\bm{\delta}^k\rangle+\frac{1}{2\eta_k}\|\bm{\delta}^k\|^2-\frac{1}{2}\|\M{V}^{k+1}-\M{W}^k\|^2_{\M{G}_k}\nonumber\\
&+\frac{1}{2}\|\M{W}-\M{W}^k\|^2_{\M{H}_k}-\frac{1}{2}\|\M{W}-\M{W}^{k+1}\|^2_{\M{H}_k}.
\end{align*}
Following a line of argument similar to the proof of \eqref{LEM22A3}, we can derive that
\begin{align}\label{LEM22A3strongly}
& \Phi(\bar{\M{X}}_K)-\Phi(\M{X}^*)+\rho\|\M{D}\bar{\M{X}}_K-\M{d}\|\nonumber\\
&\leq \frac{1}{K}\sum_{k=0}^{K-1}\big(\langle \M{X}^*-\M{X}^k,\bm{\delta}^k\rangle+\frac{\|\bm{\delta}^k\|^2}{2\eta_k}\big)+\frac{\rho^2\|\M{Q}\|}{2K}\nonumber\\
&\quad+\frac{1}{2K}\sum_{k=0}^{K-1}\frac{\|\M{X}^k-\M{X}^*\|^2-\|\M{X}^{k+1}-\M{X}^*\|^2}{\alpha_k}\nonumber\\
&\quad-\frac{1}{2K}\sum_{k=0}^{K-1}\mu\|\M{X}^k-\M{X}^*\|^2.
\end{align}
Let $\tau_k=\mu k$. Since $\frac{1}{\alpha_k}=c+\mu(k+1)=c+\tau_{k+1}$, one has
\begin{align}\label{LEM22A3strongly2}
&\sum_{k=0}^{K-1}(\frac{\|\M{X}^k-\M{X}^*\|^2-\|\M{X}^{k+1}-\M{X}^*\|^2}{\alpha_k}-\mu\|\M{X}^k-\M{X}^*\|^2)\nonumber\\
&=\sum_{k=0}^{K-1}(c+\tau_k)\|\M{X}^k-\M{X}^*\|^2-(c+\tau_{k+1})\|\M{X}^{k+1}-\M{X}^*\|^2)\nonumber\\
&\leq c\|\M{X}^0-\M{X}^*\|^2.
\end{align}
Combining \eqref{LEM22A3strongly} and \eqref{LEM22A3strongly2}, it holds that
\begin{align}\label{LEM22A3strongly3}
& \Phi(\bar{\M{X}}_K)-\Phi(\M{X}^*)+\rho\|\M{D}\bar{\M{X}}_K-\M{d}\|\nonumber\\
&\leq \frac{1}{K}\sum_{k=0}^{K-1}(\langle \M{X}^*-\M{X}^k,\bm{\delta}^k\rangle+\frac{\|\bm{\delta}^k\|^2}{2\eta_k})\nonumber\\
&\quad +\frac{c\|\M{X}^0-\M{X}^*\|^2+\rho^2\|\M{Q}\|}{2K}.
\end{align}
Since $\bm{\delta}^k=\M{g}^k-\nabla F(\M{X}^k)$, it gives that
$
\mathbb{E}(\bm{\delta}^k)=0, ~\mathbb{E}(\|\bm{\delta}^k\|^2)\leq \sigma^2,\forall k\geq0.
$
Thus, together with $\eta_k=\mu(k+1)$, it holds that
$$
\frac{1}{K}\sum_{k=0}^{K-1}\frac{\|\bm{\delta}^k\|^2}{2\eta_k}\leq\frac{\sigma^2}{2\mu K}\sum_{k=0}^{K-1}\frac{1}{k+1}\leq\frac{\sigma^2\ln K}{2\mu K}.
$$
Taking the expectation of both sides in \eqref{LEM22A3strongly3}, we get \eqref{SRATE4} immediately.
\hfill{$\blacksquare$}

\section*{Appendix F: Proof of Theorem \ref{TH3}}\label{APTH3}
{\it Proof}:
Let $v^k=(\alpha \M{D}\tr\M{\Lambda}^k+\bar{\M{X}}^k,\M{\Lambda}^k)$, $u^k=(\bar{\M{X}}^k,2\M{\Lambda}^k-\M{\Lambda}^{k-1})$, $z^k=(\M{X}^k,\M{\Lambda}^k)$, and $\bm{P}: \mathcal{W}\rightarrow \mathcal{W}:(\M{X},\M{\Lambda})\mapsto (\M{X},\alpha\M{J}\M{\Lambda})$. It follows from the iteration of the unified framework \eqref{BAS3}--\eqref{BAS2} that
\begin{subequations}\label{EQ:DYS1}
\begin{align}
z^{k+1}&=v^k-\alpha\bm{P}^{-1}\bm{B}z^{k+1},\\
u^{k+1}&=v^k-\alpha\bm{P}^{-1}((2\bm{B}+\bm{C})z^{k+1}+\bm{A}u^{k+1}),\\
v^{k+1}&=v^{k}-\alpha\bm{P}^{-1}((\bm{B}+\bm{C})z^{k+1}+\bm{A}u^{k+1}),\label{EQ:DYS1.3}
\end{align}
\end{subequations}
where $\bm{C}z^{k+1}=(\M{g}^{k+1},0)$. Let $v^*=(\alpha \M{D}\tr\M{\Lambda}^*+{\M{X}}^*,\M{\Lambda}^*)$, $u^*=({\M{X}}^*,2\M{\Lambda}^*-\M{\Lambda}^*)$, $z^*=(\M{X}^*,\M{\Lambda}^*)$, where $\M{X}^*\in\mathcal{X}^*$, and $\M{\Lambda}^*\in\mathcal{Y}^*$. It can be verified that $v^*=v^*+u^*-z^*$ and
\begin{subequations}\label{EQ:DYS2}
\begin{align}
&z^*=v^*-\alpha\bm{P}^{-1}\bm{B}z^*,\\
&u^*=v^*-\alpha\bm{P}^{-1}((2\bm{B}+\bm{C})z^*+\bm{A}u^*),\\
&v^*=v^*-\alpha\bm{P}^{-1}((\bm{B}+\bm{C})z^*+\bm{A}u^*).\label{EQ:DYS2.3}
\end{align}
\end{subequations}
From \eqref{EQ:DYS1.3} and \eqref{EQ:DYS2.3}, one can deduce that
\begin{align*}
&\|v^{k+1}-v^*\|_{\bm{P}}^2=\|v^k-v^*\|_{\bm{P}}^2-2\alpha\langle v^k-v^*,\bm{A}u^{k+1}-\bm{A}u^* \\
&+(\bm{B+C})z^{k+1}-(\bm{B+C})z^*\rangle\\
&+\alpha^2\|\bm{P}^{-1}((\bm{B+C})z^{k+1}-(\bm{B+C})z^*+\bm{A}u^{k+1}-\bm{A}u^*)\|^2_{\bm{P}}.
\end{align*}
Then, expanding the last squared $\bm{P}$-norm and using \eqref{EQ:DYS1} and \eqref{EQ:DYS2} in the inner product, it gives that
\begin{align*}
&\|v^{k+1}-v^*\|^2_{\bm{P}}=\|v^k-v^*\|_{\bm{P}}^2\\
&-2\alpha\langle z^{k+1}-z^*,(\bm{B+C})z^{k+1}-(\bm{B+C})z^*\rangle\\
&-2\alpha\langle\bm{A}u^{k+1}-\bm{A}u^*,u^{k+1}-u^*\rangle\\
&-2\alpha\langle (\bm{B+C})z^{k+1}-(\bm{B+C})z^*, \alpha \bm{B}z^{k+1}-\alpha\bm{B}z^* \rangle\\
&-2\alpha \langle \bm{A}u^{k+1}-\bm{A}u^*,\alpha((2\bm{B}+\bm{C})z^{k+1}+\bm{A}u^{k+1})\\
&\quad \quad \quad- \alpha((2\bm{B}+\bm{C})z^*+\bm{A}u^*)\rangle\\
&+\alpha^2\|\bm{P}^{-1}(\bm{B}z^{k+1}+\bm{A}u^{k+1}-\bm{B}z^*-\bm{A}u^*)\|^2_{\bm{P}}\\
&+\alpha^2\|\bm{P}^{-1}(\bm{C}z^{k+1}-\bm{C}z^*)\|^2_{\bm{P}}\\
&+2\alpha^2\langle\bm{B}z^{k+1}+\bm{A}u^{k+1}-\bm{B}z^*-\bm{A}u^*,\bm{C}z^{k+1}-\bm{C}z^*\rangle\\
&=\|v^k-v^*\|_{\bm{P}}^2-2\alpha\langle \bm{B}z^{k+1}-\bm{B}z^*,z^{k+1}-z^* \rangle\\
&-2\alpha\langle \bm{A}u^{k+1}-\bm{A}u^*,u^{k+1}-u^*\rangle\\
&-2\alpha\langle\bm{C}z^{k+1}-\bm{C}z^*,z^{k+1}-z^*\rangle\\
&+\alpha^2\|\bm{P}^{-1}(\bm{C}z^{k+1}-\bm{C}z^*)\|^2_{\bm{P}}\\
&-\alpha^2\|\bm{P}^{-1}(\bm{B}z^{k+1}+\bm{A}u^{k+1}-\bm{B}z^*-\bm{A}u^*)\|^2_{\bm{P}}.
\end{align*}
Since $\bm{A}u^{k}=(\partial R(\bar{\M{X}}^{k}),0)$, $\bm{B}z^{k}=(\M{D}\tr\M{\Lambda}^{k},\M{d}-\M{D}\M{X}^{k})$, $\bm{C}z^{k}=(\M{g}^k,0)$, $\bm{A}u^*=(\partial R({\M{X}}^*),0)$, $\bm{B}z^*=(\M{D}\tr\M{\Lambda}^*,\M{d}-\M{D}\M{X}^*)$, and $\bm{C}z^*=(\nabla F(\M{X}^*),0)$, there exist $\bm{\vartheta}^{k+1}\in\partial R(\bar{\M{X}}^{k+1})$ and $\bm{\vartheta}^*\in\partial R({\M{X}}^*)$, such that
\begin{align*}
&\|v^{k+1}-v^*\|^2_{\bm{P}}\leq \|v^k-v^*\|^2_{\bm{P}}+\alpha^2\|\M{g}^{k+1}-\nabla F(\M{X}^*)\|^2\\
&-2\alpha\langle \M{g}^{k+1}-\nabla F(\M{X}^*),\M{X}^{k+1}-\M{X}^* \rangle\\
&-2\alpha\langle \bm{\vartheta}^{k+1}- \bm{\vartheta}^*,\bar{\M{X}}^{k+1}-\bar{\M{X}}^*\rangle.
\end{align*}
Taking the expectation of both sides in the above inequality and using Assumption \ref{ASS:VR}, it holds that
\begin{align}\label{PROOF:EQ:V1}
&\mathbb{E}\|v^{k+1}-v^*\|^2_{\bm{P}}\leq \|v^k-v^*\|^2_{\bm{P}}\nonumber\\
&-2\alpha\langle \nabla F(\M{X}^{k+1})-\nabla F(\M{X}^*),\M{X}^{k+1}-\M{X}^* \rangle\nonumber\\
&+\alpha^2(2c_1D_F(\M{X}^{k+1},\M{X}^*)+c_2\sigma^2_{k+1})\nonumber\\
&-2\alpha\mathbb{E}\langle \bm{\vartheta}^{k+1}- \bm{\vartheta}^*,\bar{\M{X}}^{k+1}-\bar{\M{X}}^*\rangle.
\end{align}
Note that for any $\M{X}_1$, $\M{X}_2$
\begin{align*}
&\langle \nabla F(\M{X}_1)-\nabla F(\M{X}_2),\M{X}_1-\M{X}_2 \rangle-D_F(\M{X}_1,\M{X}_2)\\
&=F(\M{X}_2)-F(\M{X}_1)-\langle \nabla F(\M{X}_1),\M{X}_2-\M{X}_1 \rangle\geq0.
\end{align*}
It follows from \eqref{PROOF:EQ:V1} that
\begin{align*}
&\mathbb{E}\|v^{k+1}-v^*\|^2_{\bm{P}}\leq \|v^k-v^*\|^2_{\bm{P}}-2\alpha D_F(\M{X}^{k+1},\M{X}^*)\nonumber\\
&+\alpha^2(2c_1D_F(\M{X}^{k+1},\M{X}^*)+c_2\sigma^2_{k+1})\nonumber\\
&-2\alpha\mathbb{E}\langle \bm{\vartheta}^{k+1}- \bm{\vartheta}^*,\bar{\M{X}}^{k+1}-\bar{\M{X}}^*\rangle.
\end{align*}
Then, by Assumption \ref{ASS:VR}, one has
\begin{align*}
&\mathbb{E}\|v^{k+1}-v^*\|^2_{\bm{P}}+\kappa\alpha^2\mathbb{E}\sigma_{k+2}^2\\
&\leq \|v^k-v^*\|^2_{\bm{P}}+\kappa\alpha^2(1-c_3+\frac{c_2}{\kappa})\sigma^2_{k+1}\\
&-2\alpha(1-\alpha(c_1+\kappa c_4))D_F(\M{X}^{k+1},\M{X}^*)\\
&-2\alpha\mathbb{E}\langle \bm{\vartheta}^{k+1}- \bm{\vartheta}^*,\bar{\M{X}}^{k+1}-\bar{\M{X}}^*\rangle.
\end{align*}
It can be deduced that
\begin{align*}
&\mathbb{E}\|v^{k+1}-v^*\|^2_{\bm{P}}+\kappa\alpha^2\mathbb{E}\sigma_{k+2}^2\\
&\leq \|v^k-v^*\|^2_{\bm{P}}+\kappa\alpha^2(1-c_3+\frac{c_2}{\kappa})\sigma^2_{k+1}\\
&-2\alpha(1-\alpha(c_1+\kappa c_4))(D_F(\M{X}^{k+1},\M{X}^*)+\mathbb{E}D_R(\bar{\M{X}}^{k+1},\M{X}^*))
\end{align*}
Since $1-c_3+\frac{c_2}{\kappa}=1$, $\alpha\leq\frac{1}{2(c_1+\kappa c_4)}$. Let $V^k=\|v^k-v^*\|^2_{\bm{P}}+\kappa \alpha^2\sigma_{k+1}^2$. It holds that
$
\mathbb{E}V^{k+1}\leq \mathbb{E}V^k-\alpha\mathbb{E}((D_F(\M{X}^{k+1},\M{X}^*))+D_{R}(\bar{\M{X}}^{k+1},\M{X}^*)).
$
Summing it over $k=0,1,\cdots,K-1$, it gives that
$\alpha\sum_{k=0}^{K-1}\mathbb{E}((D_F(\M{X}^{k+1},\M{X}^*))+D_{R}(\bar{\M{X}}^{k+1},\M{X}^*))\leq \mathbb{E}V^0.$
Therefore, it follows from the convexity of the Bregman divergence that \eqref{VR:RATE} holds.
\hfill{$\blacksquare$}
%%%%%%%%%%%%%%%%%%%%%%%%%%%%%%%%%%%%%%%%%%%%%%%%%%%%%%%%%%%%%%%%%%%%%%%%%%%%

%%%%%%%%%%%%%%%%%%%%%%%%%%%%%%%%%%%%%%%%%%%%%%%%%%%%%%%%%%%%%%%%%%%%%%%%%%%%%%%%
\bibliographystyle{dcu}
\bibliography{../../../../bib/refference0}
%\bibliography{abbr_bibli}

\end{document}